\documentclass[12pt]{article}
\usepackage{array}
\usepackage{dsfont}
\usepackage{hyperref}
\usepackage[affil-sl]{authblk}
\usepackage{breakurl}
\usepackage{fullpage}

\usepackage{titlesec}
\titleformat{\section}
  {\Large\center\bfseries\scshape}
  {\thesection.}{.7em}{}
\titlespacing*{\section}{0pt}{3.5ex plus 0ex minus 0ex}{1.5ex plus 0ex}

\titleformat{\subsection}
  {\center\bfseries\scshape}
  {\thesubsection.}{.7em}{}
\titlespacing*{\subsection}{0pt}{3.5ex plus 0ex minus 0ex}{1.5ex plus 0ex}

\makeatletter
\g@addto@macro\normalsize{%
  \setlength\abovedisplayskip{7pt}
  \setlength\belowdisplayskip{7pt}
  \setlength\abovedisplayshortskip{7pt}
  \setlength\belowdisplayshortskip{7pt}
}
\makeatother
\usepackage[english]{babel}
\usepackage{amsfonts}
\usepackage{mathrsfs}
\usepackage{bbm}
\usepackage{amssymb}
\usepackage{amsmath}
\usepackage{todonotes}
\usepackage{enumitem}
\setlist{nolistsep}
\usepackage{amsthm}
\usepackage{amscd}
\usepackage[capitalize]{cleveref}

\newtheoremstyle{plain}
{3mm}	% Space above
{3mm}	% Space below
{\slshape}	% Body font
{}	% Indent amount
{\color{Blue}\bfseries}	% Theorem head font
{.}	% Punctuation after theorem head
{.5em}	% Space after theorem head
{}	% Theorem head spec (empty=`normal')
\newtheoremstyle{definition}
{2mm}
{2mm}
{}
{}
{\color{Blue}\bfseries}
{.}
{.5em}
{}
\theoremstyle{plain}
	
\newtheorem{Theorem}{Theorem}[section]

\newtheorem{Lemma}[Theorem]{Lemma}
\newtheorem{Proposition}[Theorem]{Proposition}
\newtheorem{Corollary}[Theorem]{Corollary}
\theoremstyle{definition}
\newtheorem{Definition}[Theorem]{Definition}
\newtheorem{Remark}[Theorem]{Remark}
\newtheorem{Example}[Theorem]{Example}

\theoremstyle{plain}
\newcounter{MainTheoremCounter}

\theoremstyle{plain}
\newtheorem*{namedthm}{\namedthmname}
\newcounter{namedthm}
\makeatletter
	\newenvironment{named}[2]
	{\def\namedthmname{#1}
	\refstepcounter{namedthm}
	\namedthm[#2]\def\@currentlabel{#1}}
	{\endnamedthm}
\makeatother
\numberwithin{equation}{section}
\usepackage{xcolor}

\definecolor{Scarlet}{rgb}{0.71, 0.11, 0.0}
\definecolor{Blue}{rgb}{0.0, 0.0, 0.0}
\definecolor{Green}{rgb}{0.39, 0.71 ,0.0}
\definecolor{Blueviolet}{rgb}{0.78, 0.11, 0.0}

\hypersetup{citecolor = Blueviolet,colorlinks,
			linkcolor = black,
			urlcolor = Blueviolet}
\newcommand{\N}{\mathbb{N}}
\newcommand{\Z}{\mathbb{Z}}
\newcommand{\R}{\mathbb{R}}
\newcommand{\C}{\mathbb{C}}
\newcommand{\Q}{\mathbb{Q}}
\newcommand{\T}{\mathbb{T}}
\newcommand{\define}[1]{\textit{#1}}

\newcommand\pile[2]{\genfrac{}{}{0pt}{}{#1}{#2}}

\renewcommand{\epsilon}{\varepsilon}
\renewcommand{\leq}{\leqslant}
\renewcommand{\geq}{\geqslant}

\renewcommand{\setminus}{\backslash}
\renewcommand{\tilde}[1]{\widetilde{#1}}
\renewcommand{\hat}[1]{\widehat{#1}}

\newcommand{\E}{\mathbb{E}}
\newcommand{\B}{\mathcal{B}}

\newcommand{\xbm}{(X,\mathcal{B},\mu)}
\newcommand{\xbmt}{(X,\mathcal{B},\mu,T)}
\renewcommand{\P}{\mathbb{P}}

\newcommand{\cZ}{\mathcal{Z}}
\renewcommand{\d}{~\mathtt{d}}

\begin{document}
% ================================================================
\allowdisplaybreaks

% =================================================TITLE&ABSTRACT
\title{%Spectral properties of multi-correlation sequences,Spectral aspects of multiple ergodic theorems
A spectral refinement of the Bergelson-Host-Kra decomposition and new multiple ergodic theorems}
\author{Joel Moreira\\
  {%
    \small
    \href{mailto:joel.moreira@northwestern.edu}{\url{joel.moreira@northwestern.edu}}
  }}
\affil{\small Department of Mathematics\\
  Northwestern University\\
  Evanston, Illinois
  }
  \author{Florian Karl Richter\\
  {%
    \small
    \href{mailto:richter.109@osu.edu}{\url{richter.109@osu.edu}}
  }}
\affil{\small Department of Mathematics\\
  The Ohio State University\\
  Columbus, Ohio
  }
\date{}
\maketitle
\begin{abstract}
We investigate how spectral properties of a measure preserving system $\xbmt$ are reflected in the multiple ergodic averages arising from that system.
For certain sequences $a:\N\to\N$ we provide natural conditions on the spectrum $\sigma(T)$ such that for all $f_1,\ldots,f_k\in L^\infty$,
\begin{equation*}
\lim_{N\to\infty}\frac{1}{N}\sum_{n=1}^N \prod_{j=1}^k T^{ja(n)}f_j=\lim_{N\rightarrow\infty}\frac{1}{N}\sum_{n=1}^N\prod_{j=1}^k T^{jn}f_j
\end{equation*}
in $L^2$-norm.
In particular, our results apply to infinite arithmetic progressions $a(n)=qn+r$, Beatty sequences $a(n)=\lfloor \theta n+\gamma\rfloor$, the sequence of squarefree numbers $a(n)=q_n$, and the sequence of prime numbers $a(n)=p_n$.

We also obtain a new refinement of Szemer{\'e}di's theorem via Furstenberg's correspondence principle.
\end{abstract}
%\tableofcontents
% ================================================================

% =========================================================SECTION
\section{Introduction}
\label{sec:intro}
% ================================================================

Let $(X,\B,\mu)$ be a probability space and let $T:X\to X$ be a measure preserving transformation.
The \define{discrete spectrum} $\sigma(T)$ of the measure preserving system
$(X,\B,\mu,T)$ is the set of eigenvalues $\theta\in\T:=\R/\Z$ for which there exists a non-zero eigenfunction
$f\in L^2(X)$ satisfying $Tf:=f\circ T=e(\theta)f$, where $e(\theta):=e^{2\pi i\theta}$.
It follows from the spectral theorem that given any functions $f,g\in L^2(X)$,
there exists a complex measure $\nu$
on the torus $\T$ such that
\begin{equation*}\label{eq_introfirst}
  \alpha(n):=\int_Xf\cdot T^ng\d\mu=\int_\T e(nx)\d\nu(x).
\end{equation*}
Decomposing the measure $\nu$ into its discrete and continuous components $\nu=\nu_d+\nu_c$, one can then represent the single-correlation sequence $\alpha(n)$ as
\begin{equation}
\label{eqn:herglotz}
\alpha(n)=\phi(n)+\omega(n),
\end{equation}
where $\phi(n)=\hat{\nu_d}(n)=\int_\T e(nx)\d\nu_d(x)$ is a \define{Bohr almost periodic sequence}\footnote{A sequence $\phi$ is \define{Bohr almost periodic} if there exists a compact abelian group $G$ (written multiplicatively), elements $a,x\in G$ and a continuous function $F\in C(G)$ such that $\phi(n)=F(a^nx)$.}
and $\omega(n)=\hat\nu_c(n)=\int_\T e(nx)\d\nu_c(x)$ is a \define{null-sequence}\footnote{A sequence $\omega$ is a \define{null-sequence} if $\lim_{N-M\to\infty}\frac1{N-M}\sum_{n=M}^N|\omega(n)|=0$.}.
It is a further consequence of the spectral theorem that for any frequency $\theta\in\T$ which does not belong to the discrete spectrum $\sigma(T)$, one has
$
\lim_{N\to\infty}\frac1N\sum_{n=1}^N\phi(n)e(-\theta n)=0
$
and hence
\begin{equation*}
\label{eqn:containment-of-spectra-2}
\lim_{N\to\infty}\frac1N\sum_{n=1}^N\alpha(n)e(-\theta n)=0.
\end{equation*}
This observation can be reformulated as
\begin{equation}
\label{eqn:containment-of-spectra}
\sigma(\alpha)\subset \sigma(T),
\end{equation}
where $\sigma(\alpha)$ denotes the \define{spectrum} of the sequence $\alpha$ in the sense introduced by Rauzy in \cite{Rauzy76}, which we recall now.

\begin{Definition}\label{def_spectrum}
  The \define{spectrum} $\sigma(\eta)$ of an arbitrary bounded sequence $\eta:\N\to\C$ is the set of frequencies $\theta\in\T$ for which
  $$\limsup_{N\to\infty}\left|\frac1N\sum_{n=1}^N\eta(n)e(-\theta n)\right|>0.$$
  Throughout this paper we are only concerned with sequences $\eta$ for which the limit supremum in the above expression is an actual limit.
\end{Definition}

Formula \eqref{eqn:containment-of-spectra} serves as the premise of our paper. Its importance lies in the many variations of the mean ergodic theorem that one can derive from it.
For instance, for any $q\in\N$, if the discrete spectrum of an ergodic system $(X,\B,\mu,T)$ is disjoint from $\{\tfrac1q,\tfrac2q,\cdots,\tfrac{q-1}q\}$ then
it follows from \eqref{eqn:containment-of-spectra} that for any $r\in\N$ and any $f\in L^2(X)$,
\begin{equation}\label{eq_toterg1}\lim_{N\to\infty}\frac1N\sum_{n=1}^N T^{qn+r}f=\int_X f \d\mu \qquad\text{in }L^2(X).
\end{equation}
More generally, for any real number $\theta>0$, if the discrete spectrum of an ergodic system $(X,\B,\mu,T)$ is disjoint from $\{\tfrac n\theta\bmod1:n\in\Z\}\setminus\{0\}$ then for any $\gamma\in\R$ and any $f\in L^2(X)$,
\begin{equation}\label{eq_totergbeatty}\lim_{N\to\infty}\frac1N\sum_{n=1}^N T^{\lfloor\theta n+\gamma\rfloor}f=\int_X f \d\mu \qquad\text{in }L^2(X).
\end{equation}
Also, invoking classical equidistribution results of Vinogradov \cite{Vinogradov47}, one can derive from \eqref{eqn:containment-of-spectra} the following ergodic theorem for totally ergodic systems: for all $f\in L^2(X)$,
\begin{equation}\label{eq_toterg2}
\lim_{N\to\infty}\frac1{\pi(N)}\sum_{\pile{p\leq N}{p~\text{prime}}}T^{p}f=\int_X f \d\mu\qquad\text{in }L^2(X),
\end{equation}
where $\pi(N)$ denotes the prime-counting function.

In this paper we seek to extend \eqref{eqn:containment-of-spectra} from single-correlation sequences to multi-correlation sequences, i.e., sequences of the form
$$\alpha(n)=\int_Xf_0\cdot T^nf_1\cdot\ldots\cdot T^{kn}f_k\d\mu$$
where $f_0,\dots,f_k\in L^\infty(X)$.
Among other things, this will allow us to derive generalizations of \eqref{eq_toterg1}, \eqref{eq_totergbeatty} and \eqref{eq_toterg2}.

The theory of multi-correlation sequences was pioneered by Furstenberg in connection with his ergodic-theoretic proof of Szemer\'edi's theorem \cite{Furstenberg77}.
A result of Bergelson, Host and Kra \cite{BHK05} offers a decomposition of multi-correlation sequences in analogy with \eqref{eqn:herglotz}:
\begin{equation}\label{eq_tremendous}
\alpha(n)=
%\int_Xf_0\cdot T^nf_1\cdots T^{kn}f_k\d\mu=
\phi(n)+\omega(n)
\end{equation}
where $\omega$ is a null-sequence and $\phi$ is a \define{nilsequence}\footnote{A sequence $\phi$ is a \define{nilsequence} if it can be approximated in $\ell^\infty$ by sequences of the form $n\mapsto F(a^nx)$, where $F$ is a continuous function on the compact homogenous space $X$ of a nilpotent Lie group $G$, $a\in G$ and $x\in X$.}.
%The Bergelson-Host-Kra decomposition allows for a deeper understanding of multi-correlation sequences by studying their nilsequence component using \define{higher order Fourier analysis} (see \cite{Tao12}), parallel to the way in which classical Fourier analysis can be used to study the Bohr almost periodic component of single-correlation sequences.
By examining the spectrum of the nilsystem from which the
nilsequence in \eqref{eq_tremendous} arises, we show that the spectrum of the multicorrelation sequence $\alpha(n)$ is contained in the discrete spectrum of its originating system $(X,\B,\mu,T)$ (cf. \cref{thm:MAIN} below).
From this we derive several multiple ergodic theorems (see
Theorems
\ref{thm_sharap},
\ref{thm_besicovitch},
\ref{thm:along-primes} and
\ref{thm_beattyprimes} in \cref{sec_results}).
As corollaries of Theorems \ref{thm_sharap} and \ref{thm:along-primes}, we obtain the following generalizations of \eqref{eq_totergbeatty} and \eqref{eq_toterg2}:
For any real number $\theta>0$, if the discrete spectrum of an ergodic system $(X,\B,\mu,T)$ is disjoint from $\{\tfrac n\theta\bmod1:n\in\Z\}\setminus\{0\}$ then for any $k\in\N$, $\gamma\in\R$ and any $f_1,\dots,f_k\in L^\infty(X)$
$$\lim_{N\to\infty}\frac1N\sum_{n=1}^N \prod_{j=1}^kT^{j\lfloor n\theta+\gamma\rfloor}f_j=\lim_{N\to\infty}\frac1N\sum_{n=1}^N \prod_{j=1}^kT^{jn}f_j \quad\text{in }L^2(X).$$
For any totally ergodic system $\xbmt$, any $k\in\N$ and any $f_1,\dots,f_k\in L^\infty(X)$,
$$
\lim_{N\to\infty}\frac1{\pi(N)}\sum_{\pile{p\leq N}{p~\text{prime}}} \prod_{j=1}^kT^{jp}f_j=\lim_{N\to\infty}\frac1N\sum_{n=1}^N \prod_{j=1}^kT^{jn}f_j \qquad\text{in }L^2(X).
$$
\paragraph{\textbf{Structure of the paper:}}
In \cref{sec_results} we formulate the main results of this paper and present relevant examples as well as applications to combinatorics. In \cref{sec:prelims} we provide the necessary background on the theory of nilmanifolds and nilsystems, which is used in the rest of the paper.

Our main technical result, \cref{thm:MAIN}, is proven in three steps: in \cref{sec_proofs}, we reduce it to the special case of nilsystems. In \cref{sec:special-case} we derive a proof of this special case conditionally on a result involving the spectrum of nilsystems.
The proof of the latter is provided in \cref{sec_main-tehnical-result}.
In the remainder of \cref{sec_proofs} and also in \cref{sec_primes} we deduce the other results stated in \cref{sec_results} from \cref{thm:MAIN}.

\paragraph{\textbf{Acknowledgements:}}
We thank V. Bergelson for many inspiring conversations, N. Frantzikinakis, B. Kra and A. Leibman for providing helpful references, and D. Glasscock for helpful comments regarding an earlier draft of this paper.
The authors also thank the referee for several pertinent suggestions which greatly improved the final version of this paper.

\section{Statement of results}
\label{sec_results}

In this section we state the main results of the paper; the proofs are presented in Sections \ref{sec_proofs} and \ref{sec_primes}.
The following theorem is our main technical result.
For a definition of nilsystems, see \cref{sec:prelims}.

\begin{Theorem}
\label{thm:MAIN}
Let $k\in\N$, let $\xbmt$ be an ergodic
measure preserving system and let $f_0,f_1,\ldots,f_k \in L^\infty(X)$.
For every $\epsilon>0$ there exists a decomposition of the form
    $$
\alpha(n):=\int f_0\cdot T^{n}f_1\cdot\ldots\cdot T^{kn}f_k\d\mu = \phi(n)+\omega(n)+\gamma(n),
$$
where $\omega(n)$ is a null-sequence, $\gamma$ satisfies $\|\gamma\|_\infty<\epsilon$ and $\phi(n)=F(R^ny)$ for some $F\in C(Y)$ and $y\in Y$, where $(Y,R)$ is a $k$-step nilsystem whose discrete spectrum is contained in the discrete spectrum of $\xbmt$.
\end{Theorem}

\begin{Remark}
A natural question is whether analogues of \cref{thm:MAIN} hold for commuting transformations or for polynomial iterates. However these extensions seem to be out of reach by the methods used in the current paper.
\end{Remark}

The following is an immediate corollary of \cref{thm:MAIN} and generalizes \eqref{eqn:containment-of-spectra}.
\begin{Corollary}
  Under the same assumptions as \cref{thm:MAIN}, the spectrum $\sigma(\alpha)$ of the multi-correlation sequence $\alpha$ (see \cref{def_spectrum}) is contained in the discrete spectrum $\sigma(T)$ of the system $\xbmt$.
\end{Corollary}

From \cref{thm:MAIN} we derive various multiple ergodic theorems.
The first theorem we derive this way is an extension of \eqref{eq_toterg1}.
An equivalent result was proven by Frantzikinakis in \cite{Frantzikinakis04}.
In the following we will use $\langle I\rangle$ to denote the subgroup of $\T$
generated by a subset $I\subset\T$.
Subsets of $\R$ are tacitly identified with their projections $\bmod~1$ onto $\T$.

\begin{Theorem}[cf. {\cite[Theorem 6.4]{Frantzikinakis04}}]
\label{thm:total-ergodicity-of-all-orders}
 Let $q,r\in\N$, and let $\xbmt$ be an ergodic measure preserving system
 whose discrete spectrum $\sigma(T)$ satisfies
 $\sigma(T)\cap\left\langle q^{-1}\right\rangle=\{0\}$.
 For any $f_1,\ldots,f_k \in L^\infty(X)$,
\begin{equation}
\label{eqn:q2-HK}
\lim_{N-M\rightarrow\infty}\frac{1}{N-M}\sum_{n=M}^N \prod_{j=1}^k T^{j(qn+r)}f_j
=
\lim_{N\rightarrow\infty}\frac{1}{N}\sum_{n=1}^N\prod_{j=1}^k T^{jn}f_j,
\end{equation}
where convergence takes place in $L^2(X)$.
In particular, if $\xbmt$ is totally ergodic, then equation \eqref{eqn:q2-HK}
holds for all $q,r\in \N$.
\end{Theorem}

The case $k=3$ of \cref{thm:total-ergodicity-of-all-orders}
was proven by Host and Kra
in \cite{HK02}. In the same paper
\cref{thm:total-ergodicity-of-all-orders} for $k>3$ was posed as a question
(\cite[Question 2]{HK02}).

\cref{thm:total-ergodicity-of-all-orders} features multi-correlation sequences along infinite arithmetic progressions $q\Z+r$.
The next theorem is a generalization in which
infinite arithmetic progressions are replaced by more general Beatty sequences $\{\lfloor \theta n+\gamma\rfloor:n\in\Z\}$.

\begin{Theorem}\label{thm_sharap}
 Let $\theta,\gamma\in\R$ with $\theta>0$, and let $\xbmt$ be an ergodic measure preserving system whose discrete spectrum $\sigma(T)$ satisfies $\sigma(T)\cap\left\langle\theta^{-1}\right\rangle=\{0\}$. For any $f_1,\ldots,f_k \in L^\infty(X)$,
\begin{equation}
\label{eqn:q3-HK}
\lim_{N-M\rightarrow\infty}\frac{1}{N-M}\sum_{n=M}^N \prod_{j=1}^k T^{j\lfloor \theta n+\gamma\rfloor}f_j
=
\lim_{N\rightarrow\infty}\frac{1}{N}\sum_{n=1}^N\prod_{j=1}^k T^{jn}f_j,
\end{equation}
where convergence takes place in $L^2(X)$. In particular, since discrete spectra are always countable, we have that for any fixed system $\xbmt$ and for almost all $\theta >0$ equation \eqref{eqn:q3-HK} holds for all $\gamma\in \R$.
\end{Theorem}

\cref{thm_sharap}, together with a standard application of Furstenberg's correspondence principle (cf. \cite[Proposition 3.1]{BHK05}), implies the following combinatorial result.
Recall that the upper density of a set $A\subset\N$ is defined by
$$\bar d(A):=\limsup_{N\to\infty}\frac{\big|A\cap\{1,\dots,N\}\big|}N.$$
\begin{Corollary}\label{cor_statements}
  Let $k\in\N$, and let $A\subset\N$ have positive upper density.
  Then for almost every $\theta\in\R$ and every $\gamma\in\R$ there exists a $k$-term arithmetic progression in $A$ with common difference in the Beatty sequence $\{\lfloor \theta n+\gamma\rfloor:n\in\N\}$.
\end{Corollary}

 In fact, under the assumptions of \cref{cor_statements}, there are many arithmetic progressions contained in $A$ with common difference in a Beatty sequence.
  More precisely, there is a syndetic\footnote{A subset $S\subset\N$ is called \define{syndetic} if it has bounded gaps, more precisely if there exists $r\in\N$ such that any interval $\{n+1,\dots,n+r\}$ of length $r$ contains at least one element of $S$.} set $S\subset\N$ such that for every $n\in S$, there exist a set $A_n\subset A$ with positive upper density and with the property that for every $m\in A_n$, the set $\big\{m,m+\lfloor \theta n+\gamma\rfloor,m+2\lfloor \theta n+\gamma\rfloor,\dots,m+k\lfloor \theta n+\gamma\rfloor\big\}$ is contained in $A$.

Recall that a bounded  sequence $\phi:\N\to\C$ is \define{Besicovitch almost periodic} if for every $\epsilon>0$, there exists a trigonometric polynomial $\rho(n)=\sum_{j=1}^tc_je(\theta_jn)$, where $t\in\N$, $\theta_1,\dots,\theta_t\in\T$ and $c_1,\dots,c_t\in\C$, such that
\begin{equation}
\label{eq:bes}
\limsup_{N\to\infty}\frac1N\sum_{n=1}^N\left|\phi(n)-\rho(n)\right|<\epsilon.
\end{equation}

The indicator function $\phi(n)=1_B(n)$ of a Beatty sequence $B=\{\lfloor\theta m+\gamma\rfloor:m\in\N\}$ is a Besicovitch almost periodic sequence with spectrum contained in the subgroup $\langle \theta^{-1}\rangle$.
Thus, sacrificing the uniformity in the Ces\`aro averages on the left hand side of \eqref{eqn:q2-HK} and \eqref{eqn:q3-HK}, one can extend Theorems \ref{thm:total-ergodicity-of-all-orders} and \ref{thm_sharap}
as follows.

\begin{Theorem}\label{thm_besicovitch}
Let $\phi$ be a Besicovitch almost periodic sequence, and let $\xbmt$ be a measure preserving system whose discrete spectrum satisfies $\sigma(T)\cap\langle\sigma(\phi)\rangle=\{0\}$. For any $f_1,\ldots,f_k \in L^\infty(X)$,
$$
\lim_{N\rightarrow\infty}\frac{1}{N}\sum_{n=1}^N\phi(n)\prod_{j=1}^k T^{jn}f_j
=
\left(\lim_{N\to\infty}\frac1N\sum_{n=1}^N\phi(n)\right)\left(\lim_{N\rightarrow\infty}\frac{1}{N}\sum_{n=1}^N\prod_{j=1}^k T^{jn}f_j\right)
$$
in $L^2(X)$.
\end{Theorem}
\begin{Remark}
    \cref{thm_besicovitch} is not true for uniform Ces\`aro averages.
    One way of obtaining a version of \cref{thm_besicovitch} with uniform Ces\`aro averages is by replacing Besicovitch almost periodic sequences with Weyl almost periodic sequences\footnote{A bounded  sequence $\phi:\N\to\C$ is \define{Weyl almost periodic} if for every $\epsilon>0$, there exists a trigonometric polynomial $\rho(n)$ such that
$\limsup_{N-M\to\infty}\frac1{N-M}\sum_{n=M}^N\left|\phi(n)-\rho(n)\right|<\epsilon.$}. In fact, one can easily modify the proof of \cref{thm_besicovitch} given below to obtain a proof of this variation.
\end{Remark}
An interesting application of \cref{thm_besicovitch} concerns the sequence $(q_n)_n$ of squarefree numbers. Since the indicator function of the set of squarefree numbers is Besicovitch almost periodic with rational spectrum (cf. \cref{sec:prelim-spectrum-of-sequences}), it follows that
for any totally ergodic $\xbmt$,
  \begin{equation*}
\lim_{N\rightarrow\infty}\frac{1}{N}\sum_{n=1}^N \prod_{j=1}^k T^{jq_n}f_j
=
\lim_{N\rightarrow\infty}\frac{1}{N}\sum_{n=1}^N\prod_{j=1}^k T^{jn}f_j.
\end{equation*}

By combining \cref{thm:MAIN} with results of Green, Tao and Ziegler \cite{Green_Tao10,GT12-2,GTZ12} on the asymptotic Gowers uniformity of the von Mangoldt function, we obtain the following multiple ergodic theorem along primes for totally ergodic systems.

\begin{Theorem}
\label{thm:along-primes}
Let $k\in\N$ and let $\xbmt$ be a totally ergodic system.
For every $f_1,\ldots,f_k \in L^\infty(X)$,
\begin{equation*}
\lim_{N\rightarrow\infty}
\frac{1}{\pi(N)}\sum_{\substack{p\leq N,\\p~\text{prime}}}~
\prod_{j=1}^k T^{jp}f_j
=
\lim_{N\rightarrow\infty}\frac{1}{N}\sum_{n=1}^N~\prod_{j=1}^k T^{jn}f_j\qquad\text{in }L^2(X).
\end{equation*}
\end{Theorem}
The case $k=3$ of \cref{thm:along-primes} was obtained by Frantzikinakis, Host and Kra in \cite[Theorem 5]{Frantzikinakis_Host_Kra07}.
In the same paper they outline the proof of \cref{thm:along-primes} in full generality, conditional on the then unknown \cref{thm:W-trick}.

Using \cref{thm_besicovitch}, we obtain a strengthening of the above result involving primes in Beatty sequences. Let $\P(\theta,\gamma):=\P\cap\{\lfloor n\theta+\gamma\rfloor:n\in\N\}$ and let $\pi_{\theta,\gamma}(N):=\big|\{1,\dots,N\}\cap\P(\theta,\gamma)\big|$.

\begin{Theorem}\label{thm_beattyprimes}
  Let $\theta,\gamma\in\R$ with $\theta>0$ irrational and let $\xbmt$ be a measure preserving system whose discrete spectrum $\sigma(T)$ satisfies $\sigma(T)\cap\left\langle \Q,\theta^{-1}\right\rangle=\{0\}$.  For every $f_1,\dots,f_k\in L^\infty(X)$,
  \begin{equation}\label{eq_thm_beattyprimes}
    \lim_{N\to\infty}\frac1{\pi_{\theta,\gamma}(N)}\sum_{\substack{p\leq N,\\p\in\P(\theta,\gamma)}}~
\prod_{j=1}^k T^{jp}f_j
=
\lim_{N\rightarrow\infty}\frac{1}{N}\sum_{n=1}^N~\prod_{j=1}^k T^{jn}f_j \qquad\text{in }L^2(X).
  \end{equation}
In particular, if $\xbmt$ is totally ergodic then for almost all $\theta>0$ equation \eqref{eq_thm_beattyprimes} holds for all $\gamma\in\R$.
\end{Theorem}

% =========================================================SECTION
\section{Preliminaries}
\label{sec:prelims}
% ================================================================

In this section we give an overview of the theory of nilspaces and nilmanifolds.

% =======================================================SUBSECTION
\subsection{Nilmanifolds and sub-nilmanifolds}
\label{sec:nilmanifolds}
% ================================================================
Let $G$ be a Lie group with identity $1_G$.
The \define{lower central series} of $G$ is the sequence
\[
G=G_1 \trianglerighteq G_2\trianglerighteq G_3\trianglerighteq \ldots\trianglerighteq\{1_G\}
\]
where $G_{i+1}:=[G_i,G]$ is, as usual, the subgroup of $G$ generated by all the commutators $aba^{-1}b^{-1}$ with $a\in G_i$ and $b\in G$.
If $G_{s+1}=\{1_G\}$ for some finite $s\in\N$ we say that $G$ is \define{($s$-step) nilpotent}.
Each $G_i$ is a closed normal subgroup of $G$ (cf. \cite[Section 2.11]{Leibman05}).

Given a nilpotent Lie group $G$ and a uniform\footnote{A closed subgroup
$\Gamma$ of $G$ is called \define{uniform} if $G/\Gamma$ is compact or, equivalently,
if there exists a compact set $K$ such that $K\Gamma = G$.}
and discrete subgroup $\Gamma$ of $G$, the quotient space
$G/\Gamma$ is called a \define{nilmanifold}.
Naturally, $G$ acts continuously and transitively on $G/\Gamma$ via
left-multiplication.

Any element $g\in G$ with the property that $g^n\in \Gamma$
for some $n\in \N$ is called \define{rational} (or \define{rational with respect to $\Gamma$)}.
A closed subgroup $H$ of $G$ is then called \define{rational}
(or \define{rational with respect to $\Gamma$}) if rational elements are dense in $H$.
For example, the subgroups $G_j$ in the lower central series of $G$ are rational with respect to any uniform and discrete subgroup
$\Gamma$ of $G$. (A proof of this fact can be found in
\cite[Corollary 1 of Theorem 2.1]{Raghunathan72} for connected $G$
and in \cite[Section 2.11]{Leibman05} for the general case.)

\begin{Remark}
\label{rem:rational-subgroup-equivalences}
It is shown in \cite{Leibman06}
that a closed subgroup $H$ is rational if and only if $H\cap \Gamma$
is a uniform discrete subgroup of $H$ if and only if
$H\Gamma$ is closed in $G$.
\end{Remark}

If $X=G/\Gamma$ is a nilmanifold, then a \define{sub-nilmanifold} $Y$ of $X$ is any
closed set of the form $Y=Hx$, where $x\in X$ and where $H$ is a closed
subgroup of $G$.
It is not true that for every closed subgroup $H$ of $G$ and for
every element $x=g\Gamma$ in $X=G/\Gamma$ the set $Hx$ is a sub-nilmanifold
of $X$; as a matter of fact, from \cref{rem:rational-subgroup-equivalences}
it follows that $Hx$ is closed in $X$ (and hence a sub-nilmanifold) if and only if the subgroup
$g^{-1}Hg$ is rational with respect to $\Gamma$.

% =======================================================SUBSECTION
\subsection{Nilsystems and their dynamics}
\label{sec:nilsystems}
% ================================================================

Let $G$ be a $s$-step nilpotent Lie group and let
$X=G/\Gamma$ be a nilmanifold.
In the following we will use $R$ or ($R_a$ if we want to emphasize the dependence on $a$)
to denote the translation by a fixed element $a\in G$, i.e. $R:x\mapsto ax$.
The map $R:X\rightarrow X$ is called a
\define{nilrotation} and the pair $(X,R)$ is called a \define{($s$-step) nilsystem}.

Every nilmanifold $X=G/\Gamma$ possesses a unique
$G$-invariant probability measure called the \define{Haar measure on X}
(cf. \cite[Lemma 1.4]{Raghunathan72}).
We will use $\mu_X$ to denote this measure.

Let us state some classical results regarding the dynamics of nilrotations.
\begin{Theorem}
[see \cite{AGH63, Parry69} in the case of connected $G$ and \cite{Leibman05} in the general case]
\label{thm:dynamics-nilrotation}
Suppose $(X,R)$ is a nilsystem. Then the following are equivalent:
\begin{enumerate}	
[label=(\roman*),ref=(\roman*),leftmargin=*]
\item\label{item:nilrotation-i}
$(X,R)$ is transitive\footnote{A topological dynamical system $(X,T)$ is called
\define{transitive} if there exists at least one point with dense orbit.};
\item\label{item:nilrotation-ii}
$(X,\mu_X,R)$ is ergodic;
\item\label{item:nilrotation-iii}
$(X,R)$ is strictly ergodic\footnote{A topological dynamical system
$(X,T)$ is called \define{strictly ergodic} if there exists a unique
$T$-invariant probability measure on $X$ and additionally the orbit of
every point in $X$ is dense.};
\end{enumerate}
Moreover, the following are equivalent:
\begin{enumerate}	
[label=(\roman*),ref=(\roman*),leftmargin=*]
\setcounter{enumi}{3}
\item
$X$ is connected and $(X,\mu_X,R)$ is ergodic.
\item
$(X,\mu_X,R)$ is totally ergodic.
\end{enumerate}
\end{Theorem}

A theorem by Lesigne \cite{Lesigne91}
asserts that for any nilmanifold $X=G/\Gamma$ with connected $G$
and any $b\in G$ the closure of the set $b^\Z x:=\{b^nx:n\in\Z\}$ is a sub-nilmanifold of $X$.
(Actually, he shows that the sequence $(b^n x)_{n\in\N}$
equidistributes with respect to the Haar measure on some sub-nilmanifold of $X$, but
in virtue of \cref{thm:dynamics-nilrotation} these two assertions are equivalent.)
Leibman has extended this result as follows.
\begin{Theorem}[{\cite[Corollary 1.9]{Leibman05-3}}]
\label{thm:Leibman05-3-Cor1.9}
Let $G$ be a nilpotent Lie group and let $\Gamma\subset G$ be a uniform and discrete subgroup.
Assume $Y$ is a connected sub-nilmanifold
of $X=G/\Gamma$ and $b\in G$.
Then $\overline{b^\Z Y}:=\overline{\bigcup_{n\in\Z}b^nY}$ is a disjoint union of finitely many
connected sub-nilmanifolds of $X$.
\end{Theorem}

% ================================================================
\subsection{The Kronecker factor of a nilsystem}
\label{sec_kronecker}
% ================================================================

Let $X=G/\Gamma$ be a nilmanifold and let $L$ be a normal, closed and rational subgroup of $G$.
Since $L\Gamma$ is closed, the quotient topology on $L\backslash X\cong G/L\Gamma$ is Hausdorff and the map $\eta: X\to L\backslash X$ that sends elements $x\in X$ to their right cosets
$Lx$ is continuous and commutes with the action of $G$.
Therefore the nilsystem $(L\backslash X, R)$ is a factor of
$(X,R)$ with factor map $\eta$.

An important tool in studying equidistribution of orbits on
nilmanifolds is a theorem by Leon Green
(see \cite{AGH63,Green61,Parry70}).
In \cite{Leibman05} Leibman offers a refinement of
this classical result of Green, a special case of which we state
now. Here and throughout the text we denote by $G^\circ$ the connected component of $G$ containing the group identity $1_G$.

\begin{Theorem}[cf. {\cite[Theorem 2.17]{Leibman05}}]
\label{thm:Leibman05-Thm2.17}
Let $X=G/\Gamma$ be a connected nilmanifold, let $u\in G$ and
let $N=\big[\langle G^\circ, u\rangle,\langle G^\circ, u\rangle\big]$,
where $\langle G^\circ, u\rangle$ denotes the group
generated by $G^\circ$ and $u$.
Then $R_u$ is ergodic on $X$ if and only if
$R_u$ is ergodic on $N\backslash X$.
\end{Theorem}
Note that in \cref{thm:Leibman05-Thm2.17} it is not explicitly stated
but implied that $N$ is a normal, closed and rational subgroup of
$G$ and hence
the factor space $N\backslash X$ is well defined.

Given a measure preserving dynamical system $\xbmt$
let $\mathcal{K}$ denote the smallest sub-$\sigma$-algebra
of $\mathcal{B}$ such that any eigenfunction of $\xbmt$
becomes measurable with respect to $\mathcal{K}$.
The resulting factor system $(X,\mathcal{K},\mu,T)$ is called the
\define{(measure-theoretic) Kronecker factor} of $\xbmt$.

The following corollary of
\cref{thm:Leibman05-Thm2.17} describes the Kronecker factor
of a connected ergodic nilsystem.

\begin{Corollary}
\label{cor:Thm2.17-kronecker}
Let $X=G/\Gamma$ be a connected nilmanifold, let $u\in G$ and
assume $R_u$ is ergodic.
Define $N:=\big[\langle G^\circ, u\rangle,\langle G^\circ, u\rangle\big]$.
Then the Kronecker factor of $(X,R_u)$
is $(N\backslash X,R_u)$.
\end{Corollary}

For the proof of \cref{cor:Thm2.17-kronecker} it will be convenient
to recall the definition of  vertical characters:
Let $G/\Gamma$ be a connected nilmanifold and let
$G=G_1 \trianglerighteq G_2\trianglerighteq
\ldots\trianglerighteq G_{s}\trianglerighteq\{1_G\}$
be the lower central series of $G$.
The quotient $T:=G_s/(\Gamma\cap G_s)$ is a connected compact
abelian group and hence isomorphic to a torus $\T^d$.
We call $T$ the \define{vertical torus} of $G/\Gamma$.
Since $G_s$ is contained in the center of $G$,
the vertical torus $T$ acts naturally on $G/\Gamma$.
A measurable function $f\in L^2(G/\Gamma)$ is called
a \define{vertical character} if there exists
a continuous group character $\chi$ of $T$ such that
$f(tx)=\chi(t)f(x)$ for all $t\in T$ and almost every $x\in X$.

\begin{proof}[Proof of \cref{cor:Thm2.17-kronecker}]
Notice that the nilsystem $(X, R_u)$ is isomorphic to the nilsystem $(X',R_u)$, where $X':= \langle G^\circ, u\rangle/(\Gamma\cap \langle G^\circ, u\rangle)$. We can therefore assume without loss of generality that $G= \langle G^\circ, u\rangle$.
We proceed by induction on the nilpotency class of $G$.
Suppose $G$ is a $s$-step nilpotent Lie group.
If $s=1$, $G$ is abelian and the result is trivial.
Next, assume that $s>1$ and that
\cref{cor:Thm2.17-kronecker} has already been proven for all
nilpotent Lie groups of step $s-1$.

Observe that $N\backslash X$ is a compact group and hence $(N\backslash X,R_u)$ is contained in the Kronecker factor of $(X,R_u)$. It thus suffices to show that for all eigenfunctions $f$ of the system $(X,R_u)$ one has
\begin{equation}
\label{eqn:inv-of-eigenfunction}
\forall v\in N\qquad f\circ R_v=f\quad \text{in }L^2(X,\mu_X).
\end{equation}

Let $\theta\in\T$ be an eigenvalue of the Koopman operator associated with $R_u$, let $E_\theta\subset L^2(X,\mu_X)$ be its (non-trivial) eigenspace and let $f\in E_\theta$.

Let $T:=G_s/(\Gamma\cap G_s)$ denote the vertical torus of
$X=G/\Gamma$ and note that the action of $T$ on $X$
commutes with the action of $R_u$.
In particular, $T$ leaves
the eigenspace $E_\theta$ invariant. It thus follows from
the Peter-Weyl theorem that $E_\theta$ decomposes into
a direct sum of eigenspaces for the Koopman representation of $T$.
In other words, any $R_u$-eigenfunction $f\in E_\theta$
can be further decomposed into a sum of vertical characters that
are also contained in $E_\theta$.
It therefore suffices to establish
\eqref{eqn:inv-of-eigenfunction} in the special case where
$f\in E_\theta$ is a vertical character.

Now assume $f\in E_\theta$, $\chi$ is a group character of $T$
and $f(tx)=\chi(t)f(x)$
for all $t\in T$ and almost every $x\in X$.
We distinguish two cases; the first case where $\chi$
is trivial and the second case where
$\chi$ is non-trivial.

Let us first assume that $\chi$ is trivial, i.e. $\chi(t)=1$
for all $t\in T$. This implies that $f$ is $G_s$ invariant.
Let $G'$ denote the nilpotent Lie group $G/G_s$
and let $\xi:G\to G'$ denote the natural quotient map.
We define $\Gamma':=\xi(\Gamma)$, which is a uniform and
discrete subgroup of $G'$, and we define $X':=G'/\Gamma'$.
Since $f$ is $G_s$ invariant it can be identified with
a function $f'$ on the nilmanifold $X'$ and $f'$ is then
an eigenfunction for $R_{u'}$, where $u'=\xi(u)$.
Since $G'$ is an $(s-1)$-step nilpotent Lie group,
we can invoke the induction hypothesis and conclude that
\begin{equation}
\label{eqn:inv-of-eigenfunction-2}\forall v'\in N':=\xi(N)\qquad f'\circ R_{v'}=f'\quad \text{in }L^2(X',\mu_{X'}).
\end{equation}
However, $f$ is $G_s$ invariant, and therefore
\eqref{eqn:inv-of-eigenfunction-2}
implies \eqref{eqn:inv-of-eigenfunction}.

Now assume that $\chi$ is non-trivial. Since $T$ is connected,
any non-trivial character has full range in the unit circle. In particular, there
exists $t\in T$ such that $\chi(t)=e(-\theta)$.
Pick any element $g\in G_s$ such that $g(\Gamma\cap G_s)=t$
and define $b:= ug$.
Then from $R_u f= e(\theta) f$ and from
$R_g f=e(-\theta)f$
it follows that $R_b f =f$. Also, note that since
the actions of $R_u$ and $R_b$ on the factor
$N\backslash X$ are identical (because $G_s\subset N$),
it follows from the ergodicity of $R_u$
that $R_b$ acts ergodically on $N\backslash X$.
Finally, the groups
$N=[G,G]=[\langle G^\circ, u\rangle,\langle G^\circ, u\rangle]$
and
$[\langle G^\circ, b\rangle,\langle G^\circ, b\rangle]$
are identical and hence it follows from
\cref{thm:Leibman05-Thm2.17} that the ergodicity of $R_b$
lifts from $N\backslash X$ to $X$.
We conclude that $f$ has to be a constant function, thereby
satisfying \eqref{eqn:inv-of-eigenfunction}.
\end{proof}

\subsection{Spectrum of almost periodic sequences}
\label{sec:prelim-spectrum-of-sequences}
In this section we collect a few facts about the spectrum of almost periodic sequences; we refer the reader to the book of Besicovitch \cite{Besicovitch55} for a complete treatment on the theory of almost periodic sequences.

Almost periodic sequences were first introduced by Bohr
in \cite{Bohr25}. In his second paper on this subject \cite{Bohr25b} he proves that any
Bohr almost periodic sequence can be approximated uniformly by
trigonometric polynomials whose frequencies are all contained in the
spectrum $\sigma(\phi)$. This theorem is known as
Bohr's approximation theorem.
An analogue of Bohr's approximation theorem for Besicovitch almost periodic
sequences was later obtained by Besicovitch.
He showed that the spectrum of a Besicovitch almost periodic sequence is at most countable and then proved that any
Besicovitch almost periodic sequence can be approximated in the Besicovitch seminorm by
trigonometric polynomials whose frequencies are all contained in the
spectrum $\sigma(\phi)$.
The precise statement of Besicovitch's result is as follows\footnote{In his book, Besicovitch only deals with continuous almost periodic functions (i.e., almost periodic functions with domain $\R$), but the proof of \cite[Theorem II.8.2$^\circ$]{Besicovitch55} also works for discrete almost periodic sequences (i.e., almost periodic functions with domain $\Z$); see also \cite[Section 3]{BL85}.}.

\begin{Theorem}[cf. {\cite[Theorem II.8.2$^\circ$ (page 105)]{Besicovitch55}}]\label{lemma_besicovitchspectrum}
  Let $\phi$ be a Besicovitch almost periodic sequence with spectrum $\sigma(\phi)$.
  Then for every $\epsilon>0$ there exists a trigonometric polynomial $\rho(n)=\sum_{i=1}^tc_ie(\theta_in)$ with
$c_1,\ldots,c_t\in\C$ and
  $\theta_1,\ldots,\theta_t\in\sigma(\phi)$ such that
  \begin{equation}\label{eq_besicovitchalmostperiodic}
  \limsup_{N\to\infty}\frac1N\sum_{n=1}^N\left|\phi(n)-\rho(n)\right|<\epsilon.
  \end{equation}
\end{Theorem}

We will also make use of the following lemma regarding the spectrum of the product of two Besicovitch almost periodic sequences.

\begin{Lemma}\label{lemma_productbesicovitch}
Let $f_1,f_2:\N\to\C$ be bounded Besicovitch almost periodic sequences.
The the product $f_1\cdot f_2$ is also Besicovitch almost periodic with spectrum $\sigma(f_1\cdot f_2)\subset\langle \sigma(f_1)\cup\sigma(f_2)\rangle$.
\end{Lemma}
\begin{proof}
  In view of \cref{lemma_besicovitchspectrum}, we can approximate each $f_i$ with a trigonometric polynomial $\rho_i$ whose spectrum is contained in $\sigma(f_i)$.
  Observe that the product $\rho_1\rho_2$ is a trigonometric polynomial with spectrum contained in $\langle \sigma(f_1)\cup\sigma(f_2)\rangle$.
  Finally, it is not hard to show that $\rho_1\rho_2$ approximates $f_1f_2$, which finishes the proof.
\end{proof}

\section{Proving \cref{thm:MAIN} for the special case of nilsystems}
\label{sec:special-case}

In this section we will prove
\cref{thm:MAIN}
for the special case of nilsystems.
This will serve as an
important intermediate step in obtaining
\cref{thm:MAIN} in its full
generality.

\begin{Theorem}
\label{thm:tot-erg-nilsystems}
Let $k\in\N$, let $(X,R)$ be an ergodic $k$-step nilsystem
and let $f_0,f_1,\ldots,f_k \in C(X)$.
 Then
\begin{equation}\label{eq_thm:tot-erg-nil}
  \int f_0\cdot R^nf_1R^{2n}f_2\cdot\ldots\cdot R^{kn}f_k\d\mu_X=\phi(n)+\omega(n),
\end{equation}
where $\omega(n)$ is a null-sequence and $\phi(n)=F(S^ny)$ for some $F\in C(Y)$ and $y\in Y$, where $(Y,S)$ is a $k$-step nilsystem whose discrete spectrum $\sigma(Y,S)$ is contained in $\sigma(X,R)$, the discrete spectrum of $(X,R)$.
\end{Theorem}

The main new ingredient used in the proof of \cref{thm:tot-erg-nilsystems} is the
following result.

\begin{Theorem}
\label{cor:thm3.4}
Let $k\in\N$, let $X$ be a connected nilmanifold and let $R:X\to X$
be an ergodic nilrotation.
Define $S:= R\times R^2\times \ldots\times R^k$ and
\begin{equation*}
Y_{X^\Delta}:=\overline{\big\{S^n(x,x,\ldots, x):
x\in X,~n\in\Z\big\}}\subset X^{k}.
\end{equation*}
Then $\sigma(X,R)=\sigma(Y_{X^\Delta}, S)$.
\end{Theorem}

The proof of \cref{cor:thm3.4} is postponed to \cref{sec_main-tehnical-result}.

Most of the ideas used in the rest of the proof of \cref{thm:tot-erg-nilsystems} were already present in \cite{BHK05} and \cite{Leibman10}.
For completeness, we repeat the same arguments here, adapting them to our situation as needed.

Let $X=G/\Gamma$ be a nilmanifold and let $\pi:G\to X$ denote the natural projection of $G$ onto $X$.
We will
use $1_X$ to denote the point $\pi(1_G)$.
Consider a closed subgroup $H$ of $G$.
As noticed in \cref{rem:rational-subgroup-equivalences} the set $Y:=\pi(H)$ is a sub-nilmanifold of $X$ if and only if
$H$ is rational. Let $L$ denote the normal closure of $H$
in $G$, that is, let $L$ be the smallest normal subgroup of $G$ containing
$H$. One can show that if $H$ is closed and rational
then so is $L$ (cf. \cite{Leibman10}).
In particular, the set $Z:=\pi(L)$ is a sub-nilmanifold of
$X$ containing $Y$.
We call $Z$ the \define{normal closure} of $Y$.

Note, every sub-nilmanifold $Y=Hx$ of $X$ can be viewed as a nilmanifold on
its own and in particular it has its own Haar measure $\mu_Y$.
Moreover, for any $a\in G$, the Haar
measure of the sub-nilmanifold $aY$ coincides with the push forward of
$\mu_Y$ under $R_a$.

\begin{Proposition}[cf. {\cite[Proposition 3.1]{Leibman10}}]
\label{prop:Leibman10-Proposition-3.1}
Assume $V=H/\Gamma_H$ is a connected nilmanifold,
$\pi:H\to V$ is the natural projection of $H$ onto $V$
and $W$ is a connected sub-nilmanifold of $V$ containing
$1_V=\pi(1_H)$.
Let $b\in H$ and assume $b^\Z W$ is dense in $V$.
If $Z$ denotes the normal closure of $W$, then
for all $f\in C(V)$ we have
\begin{equation}
\label{eq:prop-3.1-eqn1}
\lim_{N-M\to\infty}\frac{1}{N-M}\sum_{n=M}^N
\left|\int_{b^n W} f\d\mu_{b^nW} -\int_{b^n Z} f\d\mu_{b^n Z}\right|=0.
\end{equation}
\end{Proposition}

\begin{Proposition}
\label{cor:Leibman10-Corollary-3.2}
Let $(X,S)$ be a nilsystem, let $W\subset X$ be a connected sub-nilmanifold containing the origin
$1_X$ and assume that $V:=\overline{S^\Z W}$ is also a connected sub-nilmanifold of $X$.
Then there exists a factor $(Y,S)$ of $(V,S)$, and a point $y\in Y$ such that for any continuous function $f\in C(X)$, there exists $F \in C(Y)$ such that
\begin{equation}
\label{eq:cor-3.2-eqn1}
\lim_{N-M\to\infty}\frac{1}{N-M}\sum_{n=M}^N
\left|\int_{S^n W} f\d\mu_{S^n W} -F(S^ny)\right|=0.
\end{equation}
\end{Proposition}

\begin{proof}
Since $V$ is invariant under $S$ and $1_X\in V$,
we can find a closed rational subgroup $H$ of $G$
such that $V=\pi(H)$.
Therefore, $\Gamma_H:=\Gamma\cap H$
is a uniform discrete subgroup of $H$ and the nilsystem
$(V,S)$ is isomorphic to $(H/\Gamma_H,S)$. In the following we will identify
$V$ with $H/\Gamma_H$ and vice versa.
Let $Z$ be the normal closure of the sub-nilmanifold $W$ in $V$
and let $L$ denote the corresponding normal subgroup of $H$ such that
$\pi(L)=Z$.

Define $Y:=L\backslash V\cong H/(L\Gamma_H)$ and let
$\eta:V\rightarrow Y$ denote the natural projection.
As explained at the beginning of \cref{sec_kronecker},
$(Y,S)$ is a well defined factor of $(V,S)$ with factor map $\eta$.

Define $y:=\eta(1_X)$ and observe that $\eta(W)=\{y\}$.
Note that for every $z=\eta(\pi(h))\in Y$, the set
$\eta^{-1}(z)=\pi(gL)$ is a sub-nilmanifold
of $Y$ and therefore it possesses a Haar measure, which we denote by
$\mu_{\eta^{-1}(z)}$.
Let $f\in C(X)$ and define the function $F$ as
\[
F(z):=\int_{\eta^{-1}(z)} f\d\mu_{\eta^{-1}(z)}.
\]
Finally, observe that $F(S^n y)= \int_{S^nZ} f\d\mu_{S^nZ}$
and so \eqref{eq:cor-3.2-eqn1} follows immediately from
\cref{eq:prop-3.1-eqn1} in
\cref{prop:Leibman10-Proposition-3.1}.
\end{proof}

To prove \cref{thm:tot-erg-nilsystems} we will also require a technical lemma:
\begin{Lemma}\label{lemma_thm:tot-erg-nilsystems}
  Let $(X,R)$ be an ergodic connected nilsystem of step $s$ and let $q\in\N$.
  Then there exists an ergodic nilsystem $(Y,S)$ of step $s$ with exactly $q$ connected components and such that the restriction of $S^q$ to each connected component of $Y$ yields a system isomorphic to $(X,R)$.
\end{Lemma}

\begin{proof}
   First, we claim that one can embed the connected nilsystem $(X,R)$ into a nilflow $(X',(R^t)_{t\in\R})$, so that $X$ is a subnilmanifold of $X'$ invariant under $R=R^1$.
  Indeed, say $X=G/\Gamma$.
  One can assume that the identity component $G^\circ$ of $G$ is simply connected, by passing to the universal cover if needed.
  Next one can use \cite[Theorem 2.20]{Raghunathan72} to find a connected simply connected nilpotent Lie group $G'$ such that $G\subset G'$ and $\Gamma$ is a uniform discrete subgroup of $G'$.
  In particular $X$ is a sub-nilmanifold of $X':=G'/\Gamma$.
Since $G'$ is connected and simply connected, for any element $a\in G'$ the associated one-parameter subgroup $(a^t)_{t\in\R}$ is well defined (cf. \cite[Subsection 2.4]{Leibman05}).
In particular, the nilrotation $R=R_a:X\to X$ can be extended to a nilflow $(R^t)_{t\in\R}$ on $X'$ by defining $R^t x:= R_{a^t}x=a^t x$ for all $x\in X'$, $t\in\R$.

Next, consider the product nilsystem $(Y',S):=(X',R^{1/q})\times(\Z/(q\Z),+1)$, so that as a nilmanifold $Y'=X'\times\{0,1,\dots,q-1\}$ and the nilrotation $S$ is defied as $S(x,r)=(R^{1/q}x,r+1\bmod q)$.
Finally, let $Y\subset Y'$ be the orbit of $X\times\{0\}$ under $S$.
Since $S^q=R^1\times Id=T\times Id$ preserves $X\times\{0\}$, we deduce that $Y$ has precisely $q$ connected components.
In fact, $S^q$ preserves each component $X\times\{r\}$, and moreover $(X\times\{r\},S^q)$ is isomorphic to $(X,R)$ via the map $x\mapsto(R^{r/q}x,r)$.
\end{proof}

We are now ready to prove \cref{thm:tot-erg-nilsystems}.
\begin{proof}[Proof of \cref{thm:tot-erg-nilsystems}]
Using invariance of the measure $\mu_X$ under $R$ yields
$$\int f_0\cdot R^nf_1R^{2n}f_2\cdot\ldots\cdot R^{kn}f_k\d\mu_X
= \int R^nf_0\cdot R^{2n}f_1\cdot\ldots\cdot R^{(k+1)n}f_k\d\mu_X.
$$
Hence, by changing $k$ to $k+1$ and renaming the functions $f_0,f_1,\ldots,f_k$ to $f_1,f_2,\ldots,f_k$, we see that in order to prove \cref{thm:tot-erg-nilsystems} it is equivalent to prove that
for any $k\in\N$, any ergodic $(k-1)$-step nilsystem $(X, R)$
and any $f_1,f_2,\ldots,f_k \in C(X)$ we have
\begin{equation}\label{eq_thm:tot-erg-nil-new}
\int R^nf_1\cdot R^{2n}f_2\cdot\ldots\cdot R^{kn}f_k\d\mu_X=\phi(n)+\omega(n),
\end{equation}
where $\omega(n)$ is a null-sequence and $\phi(n)=F(S^ny)$ is a nilsequence coming from a $(k-1)$-step nilsystem $(Y,S)$ with $\sigma(Y,S)\subset\sigma(X,R)$.

Let $\alpha(n)$ denote the sequence
$$
\alpha(n):=\int R^nf_1\cdot R^{2n}f_2\cdot\ldots\cdot R^{kn}f_k\d\mu_X.
$$
We first deal with the case when $X$ is connected.
Let $X^{\Delta}$ be the diagonal of $X^k$ and let $S=R\times R^2\times\cdots\times R^k$. Note that since $X$ is connected, the diagonal $X^{\Delta}$ is also connected. We can write $\alpha(n)$ as
$$
\alpha(n)= \int_{S^n X^{\Delta}} f_1\otimes f_2\otimes\cdots\otimes f_k\d\mu_{S^n X^{\Delta}}.
$$
It is shown in \cite[Corollary 6.5]{Leibman10b} and also in
\cite[Corollary 2.10]{Frantzikinakis08} that if $X$ is connected then
$Y_{X^\Delta}:=\overline{\{S^n X^{\Delta}:n\in\Z\}}$ is connected.
It thus follows from
\cref{thm:Leibman05-3-Cor1.9} that $Y_{X^\Delta}$
is a sub-nilmanifold of $X^k$.
Also, from \cref{cor:thm3.4} we have that $Y_{X^\Delta}$ satisfies $\sigma(Y_{X^\Delta},S)=\sigma(X,R_a)$.
This observation allows us to apply \cref{cor:Leibman10-Corollary-3.2}
with $W= X^{\Delta}$ and $V=Y_{X^\Delta}$.
Therefore we can find a factor $(Y,S)$ of $(Y_{X^\Delta},S)$, a point $y\in Y$
and a continuous function $F\in C(Y)$ such that \eqref{eq:cor-3.2-eqn1}
is satisfied.
Observe that since $(Y,S)$ is a factor of $(Y_{X^\Delta},S)$,
the discrete spectra satisfy $\sigma(Y,S)\subset\sigma(Y_{X^\Delta},S)$ and hence
$\sigma(Y,S)\subset\sigma(X,R)$.
Besides that, \eqref{eq:cor-3.2-eqn1} can be written as
$$
\lim_{N-M\to\infty}\frac{1}{N-M}\sum_{n=M}^N
\left|\alpha(n) - F(S^n y)\right|=0.$$
Therefore, setting $\phi(n)=F(S^ny)$ and $\omega(n)=\alpha(n)-\phi(n)$, we obtain a decomposition of $\alpha(n)$ satisfying \eqref{eq_thm:tot-erg-nil-new}.

Next we deal with the case when $X$ is not connected.
Since $X$ is compact, it has a finite number of connected components $X_0,\dots,X_{q-1}$.
It is not hard to see that $R$ permutates the components $X_0,\dots,X_{q-1}$ cyclically, that is, $RX_\ell=X_{\ell+1\bmod q}$.
In particular, $R^q$ preserves each $X_\ell$ and for each $n\in\Z$ the map $R^n$ is an isomorphism between the systems $(X_\ell,R^q)$ and $(X_{\ell+n\bmod q},R^q)$.
Also, observe that for each $\ell\in\{0,\dots,q-1\}$ the system $(X_\ell,R^q)$ is totally ergodic.
This shows that the function $\sum_{\ell=0}^{q-1} e\left(\frac{\ell}{q}\right)1_{X_\ell}$ is an eigenfunction for $R$ with eigenvalue $1/q$ and so $$\sigma(X,R)\cap\Q=\left\{0,\frac1q,\frac2q,\dots,\frac{q-1}q\right\}.$$
On the other hand, an irrational point $\theta\in\T$ is an eigenvalue for $R$ if and only if $q\theta$ is an eigenvalue for $R^q$; therefore we conclude that
\begin{equation}\label{eq_prooffix2}
  \sigma(X,R)=\tfrac1q\sigma(X_0,R^q)\oplus\{1/q,\dots,(q-1)/q\}.
\end{equation}

For each $r\in\{0,\dots,q-1\}$ and $i\in\{1,\dots,k\}$ let $f_{i,r}:=R^{ir}f_i$.
For $\ell\in\{0,\dots,q-1\}$ denote by
$$\alpha_{\ell,r}(n):=\int_{X_\ell}R^{qn}f_{1,r}\cdot R^{2qn}f_{2,r}\cdot\ldots\cdot R^{kqn}f_{k,r}\d\mu_{X_\ell}$$
and observe that
$$\alpha(qn+r)=\int_X\prod_{i=1}^kR^{i(qn+r)}f_i\d\mu_X=\sum_{\ell=0}^{q-1}\int_{X_\ell}\prod_{i=1}^kR^{iqn}f_{i,r}\d\mu_{X_\ell}= \sum_{\ell=0}^{q-1}\alpha_{\ell,r}(n).$$
Since $R^\ell$ is an isomorphism between $(X_0,R^q)$ and $(X_\ell,R^q)$ we have
$$\alpha_{\ell,r}(n)=\int_{X_0}\prod_{i=1}^kR^{qin}R^\ell f_{i,r}\d\mu_{X_0}=\int_{X_0}\prod_{i=1}^kR^{qin}f_{i,r,\ell}\d\mu_{X_0}$$
where $f_{i,r,\ell}=R^\ell f_{i,r}=R^{\ell+ir}f_i$.

By applying the above argument for connected nilsystems to $(X_0,R^q)$ we find a nilsystem $(Y_0,S_0)$ with spectrum $\sigma(Y_0,S_0)\subset\sigma(X_0,R^q)$, a point $y_0\in Y_0$ and, for each $r,\ell\in\{0,\dots,q-1\}$, a function $F_{\ell,r}\in C(Y_0)$ such that
$$\lim_{N-M\to\infty}\frac{1}{N-M}\sum_{n=M}^N
\left|\alpha_{\ell,r}(n) -F_{\ell,r}(S_0^ny_0)\right|=0.$$
%Furthermore, the spectrum $\sigma(Y_0,
Letting $F_r:=\sum_{\ell=0}^{q-1}F_{\ell,r}$ we deduce that
\begin{equation}\label{eq_prooffix2-102}
\lim_{N-M\to\infty}\frac{1}{N-M}\sum_{n=M}^N
\left|\alpha(qn+r) -F_r(S_0^n y_0)\right|=0.
\end{equation}
Invoking \cref{lemma_thm:tot-erg-nilsystems} we find an ergodic nilsystem $(Y,S)$ with $q$ connected components and such that $Y_0$ can be identified with one of the connected components of $Y$ in such a way that the restriction of $S^q$ to $Y_0$ is precisely $S_0$.
The same argument which led to \eqref{eq_prooffix2} implies that
$$
\sigma(Y,S)=\tfrac{1}{q}\sigma(Y_0,S_0)\oplus\left\{0,\tfrac{1}{q},\dots,\tfrac{q-1}{q}\right\}.
$$
Together with \eqref{eq_prooffix2} and the fact that $\sigma(Y_0,S_0)\subset\sigma(X_0,R^q)$, the previous equation implies that $\sigma(Y,S)\subset\sigma(X,R)$.

Each point in $Y$ can be represented uniquely as $S^ry$ for some $r\in\{0,\dots,q-1\}$ and $y\in Y_0\subset Y$.
Let $F\in C(Y)$ be defined as $F(S^ry)=F_r(y)$, define $\phi(n):=F(S^ny_0)$ and let $\omega(n):=\alpha(n)-\phi(n)$.
This way, we obtain a decomposition as in \eqref{eq_thm:tot-erg-nil}, where $\phi(n)$ is a nilsequence coming from a $(k-1)$-step nilsystem $(Y,S)$ with $\sigma(Y,S)\subset\sigma(X,R)$. It only remains to show that $\omega(n)$ is a nullsequence:

\begin{equation*}
\begin{split}
\lim_{N-M\to\infty}\frac{1}{N-M}&\sum_{n=M}^N
\left|\omega(n)\right|
\\
=&~
\frac{1}{q}\sum_{r=0}^{q-1}
\lim_{N-M\to\infty}\frac{1}{N-M}\sum_{n=M}^N
\left|\alpha(qn+r) -F(S^{qn+r}(y_0))\right|
\\
=&~
\frac{1}{q}
\sum_{r=0}^{q-1}
\lim_{N-M\to\infty}\frac{1}{N-M}\sum_{n=M}^N
\left|\alpha(qn+r) -F_r(S_0^n y_0)\right|,
\end{split}
\end{equation*}
which together with \eqref{eq_prooffix2-102} shows that $\omega(n)$ is a null-sequence.\qedhere
\end{proof}

\section{Proofs of {\cref{thm:MAIN}} and its corollaries}\label{sec_proofs}

The purpose of this section is to derive proofs of the
main theorems, namely Theorems \ref{thm:MAIN},
%\ref{thm:total-ergodicity-of-all-orders},
\ref{thm_sharap},
%\ref{thm:irrational-nilsequence-plus-nullsequence},
\ref{thm_besicovitch}.

\subsection{Host-Kra-Ziegler factors and a proof of {\cref{thm:MAIN}}}

We start with the proof of \cref{thm:MAIN}, which we recall for convenience of the reader.

\begin{named}{\cref{thm:MAIN}}{}
Let $k\in\N$, let $\xbmt$ be an ergodic measure preserving system and let $f_0,f_1,\ldots,f_k \in L^\infty(X)$.
Then for every $\epsilon>0$ we have a decomposition of the form
    $$\int f_0\cdot T^{n}f_1\cdot\ldots\cdot  T^{kn}f_k\d\mu = \phi(n)+\omega(n)+\gamma(n),$$
where $\omega(n)$ is a null-sequence, $\gamma$ satisfies $\|\gamma\|_\infty<\epsilon$ and $\phi(n)=F(R^ny)$ for some $F\in C(Y)$ and $y\in Y$, where $(Y,R)$ is a $k$-step nilsystem whose discrete spectrum is contained in the discrete spectrum of $\xbmt$.
\end{named}

In \cref{sec:special-case} we have already established this theorem under the additional assumptions that $(X,{\mathcal B},\mu,T)$ is a nilsystem and each $f_i$ is continuous.
In order to close the gap between nilsystems and general
measure preserving systems,
we rely on the theory of the Host-Kra-Ziegler factors,
which was developed by Host and Kra in \cite{HK05} and
independently by Ziegler in \cite{Ziegler07}.
Let us briefly summarize their theory; for details we refer the reader to \cite{HK05-2,HK05,Ziegler07}.

Suppose $\xbmt$ is a measure preserving system.
For $s\in\N$ the \define{$s$-th Host-Kra-Ziegler factor} of $\xbmt$,
denoted by $\cZ_s$, is a $T$-invariant sub-$\sigma$-algebra of $\B$
which serves as a characteristic factor for multiple ergodic averages.
For every $s$, the system $(X,\cZ_s,\mu,T)$ is an inverse limit
of $s$-step nilsystems, meaning that there exists a nested
sequence of $T$-invariant sub-$\sigma$-algebras
$\cZ_s^{(1)}\subset \cZ_s^{(2)}\subset \ldots \subset\cZ_s$
such that $\bigvee_{m\geq 1}\cZ_s^{(m)}=\cZ_s$
and such that for every $m$ the system $(X,\cZ_s^{(m)},\mu,T)$ is
measurably isomorphic to an $s$-step nilsystem.

As mentioned above, the factors $\cZ_s$ are
\define{characteristic factors} for multiple ergodic averages.
In particular we have the following theorem.

\begin{Theorem}[{\cite[Corollary 4.5]{BHK05}}]
\label{thm:BHK-L-structure-theorem}
For any $k\in\N$, any ergodic measure preserving system $\xbmt$ and any
$f_0,\ldots,f_k \in L^\infty(X)$,
\[
\lim_{N-M\to\infty}
\frac{1}{N-M}\sum_{n=M}^N\left|\int_X\left(\prod_{j=0}^kT^{jn}f_j-
\prod_{j=0}^k T^{jn}\E(f_j|\cZ_k)\right)\d\mu\right| =0.
\]
\end{Theorem}

\begin{proof}[Proof of \cref{thm:MAIN}]
We are given an ergodic measure preserving system $\xbmt$
as well as bounded measurable functions
$f_0,f_1,\ldots,f_k \in L^\infty(X)$ and some $\epsilon>0$;
we seek to decompose the multi-correlation sequence
\[
\alpha(n):=\int~f_0 \cdot T^{n}f_1\cdot\ldots\cdot T^{kn}f_k\d\mu
\]
into $\alpha(n)=\varphi(n)+\omega(n)+\gamma(n)$,
satisfying the stated properties.

Using \cref{thm:BHK-L-structure-theorem}, and after changing $\alpha$ by a null-sequence if necessary, we may assume that $\E(f_j|\cZ_k)=f_j$ for all $j\in\{0,1,\ldots,k\}$.
Let $\cZ_k^{(m)}$, $m\in\N$, be a nested sequence of $T$-invariant
sub-$\sigma$-algebras of $\cZ_k$
such that $\bigvee_{m\geq 1}\cZ_k^{(m)}=\cZ_k$ and such that for every
$m$ the system $(X,\cZ_k^{(m)},\mu,T)$ is measurably isomorphic to
a compact $k$-step nilsystem.
Since $(X,\cZ_k^{(m)},\mu,T)$ is a factor of $\xbmt$,
we have the inclusion of the spectra $\sigma(X,\cZ_k^{(m)},\mu,T)\subset\sigma\xbmt$.

For each $m\in\N$, let $\alpha_m:\Z\to\R$ be defined as
\[
\alpha_m(n):=\int~\E(f_0|\cZ_k^{(m)})\cdot T^{n}\E(f_1|\cZ_k^{(m)})
\cdot\ldots\cdot T^{kn}\E(f_K|\cZ_k^{(m)})\d\mu.
\]
It follows from Doob's martingale convergence theorem (cf., for
instance, \cite[Theorem 5.4.5]{Durrett10}) that $\alpha_m(n)$ converges
to $\alpha(n)$ as $m\to\infty$, uniformly in $n$.
In other words, choosing a large enough $m\in\N$ we have that $\|\alpha-\alpha_m\|_\infty<\epsilon/2$.

Next, one can approximate the functions $\E(f_i|\cZ_k^{(m)})$
(when identified with measurable functions on the respective compact
$s$-step nilsystems) by continuous functions in $\|.\|_{L^p}$-norm,
for every $p<\infty$.
More precisely, there exists an $s$-step nilsystem
$(\tilde X,\tilde\mu,\tilde T)$, which is measurably isomorphic to $(X,\cZ_k^{(m)},\mu,T)$, and continuous functions
$\tilde f_1,\dots,\tilde f_k\in C(\tilde X)$ such that
the sequence
$$
\beta(n):=
\int_{\tilde X}\tilde f_0\cdot \tilde T^{n}\tilde f_1
\cdot\ldots\cdot \tilde T^{kn}\tilde f_k\d\tilde\mu
$$
satisfies $\|\alpha_m-\beta\|_\infty <\epsilon/2$ and hence $\gamma(n):=\alpha(n)-\beta(n)$ satisfies $\|\gamma\|_\infty<\epsilon$.

Finally, it follows from
\cref{thm:tot-erg-nilsystems} that $\beta$ can be written as
\begin{equation*}
\label{eqn:pr-of-thmB-1}
\beta(n)=\phi(n) + \omega(n),
\end{equation*}
where $\omega$ is a nullsequence and $\phi(n)=F(R^ny)$ for some $F\in C(Y)$ and $y\in Y$, where $(Y,R)$ is a $k$-step nilsystem whose discrete spectrum satisfies $\sigma(Y,R)\subset\sigma(\tilde X,\tilde\mu,\tilde T)=\sigma(X,\cZ_k^{(m)},\mu,T)\subset\sigma\xbmt$, finishing the proof.
\end{proof}

% =========================================================SECTION
\subsection{Beatty sequences and a proof of \cref{thm_sharap}}
% ================================================================
For the proof of \cref{thm_sharap}, we will use the following lemma.

\begin{Lemma}\label{lem_hl1}
Let $\theta\in\T$ and let $F_1:\T\to\T$ be Riemann integrable with $\displaystyle \lim_{N\rightarrow\infty}\frac{1}{N}\sum_{n=1}^N F_1(n\theta) =0$.
Also, let $\xbmt$ be an ergodic measure preserving system whose discrete spectrum satisfies $\sigma(T)\cap\left\langle\theta\right\rangle=\{0\}$. Then
for any $f_1,\ldots,f_k \in L^\infty(X)$ we have
\begin{equation}
\label{eqn:q4-HK}
\lim_{N-M\rightarrow\infty}\frac{1}{N-M}\sum_{n=M}^N F_1(n\theta)\prod_{j=1}^k T^{jn}f_j
=
0
\end{equation}
in $L^2$.
\end{Lemma}

\begin{proof}
It follows from \cite[Theorem 2.24]{Host_Kra09} that the limit of the left hand side in
\eqref{eqn:q4-HK} exists in $L^2$. It thus suffices to show that for all $f_0\in L^\infty(X)$ the limit
\begin{equation}\label{eq_proofbesicovitchmulticorrelation-old}
\lim_{N-M\rightarrow\infty}\frac{1}{N-M}\sum_{n=M}^N F_1(n\theta)\int~f_0\cdot T^{n}f_1\cdot\ldots\cdot T^{kn}f_k\d\mu=0.
\end{equation}

In view of \cref{thm:MAIN}, the multi-correlation sequence can be decomposed as
$$\int~f_0 \cdot T^{n}f_1\cdot\ldots\cdot  T^{kn}f_k\d\mu=\psi(n)+\omega(n)+\gamma(n),$$
where $\|\gamma\|_\infty<\epsilon$, $\omega$ is a null-sequence and $\psi(n)=F_2(R_2^ny)$, where $y\in Y_2$, $F_2\in C(Y_2)$ and $(Y_2,R_2)$ is a nilsystem whose discrete spectrum satisfies $\sigma(Y_2,R_2)\subset\sigma(X,T)$. Let $\mu_{Y_2}$ denote the
Haar measure of the nilmanifold $Y_2$.

Let us use $R_1:\T\to\T$ to denote rotation by $\theta$ on $\T$ and let $Y_1$ denote the orbit closure of $0$ under $R_1$.
Note that either $Y_1=\T$, which corresponds to the case of
irrational $\theta$, or $Y_1$ is a finite subgroup of $\T$, which corresponds to the case of rational $\theta$. Either way,
$Y_1$ is a closed subgroup of $\T$ and we use $\mu_{Y_1}$
to denote its Haar measure.

Putting everything together we can now rewrite the left hand side of \eqref{eq_proofbesicovitchmulticorrelation-old} (up to a loss of $\epsilon$) as
\begin{equation*}
  \lim_{N\rightarrow\infty}\frac{1}{N}\sum_{n=1}^NF_1(R_1^n0)F_2(R_2^ny)=0.
\end{equation*}
Since the discrete spectrum $\sigma(Y_1,R_1)$ of the system $(Y_1,R_1)$
is given by the group generated by $\theta$, we deduce that the systems $(Y_1,R_1)$ and $(Y_2,R_2)$ have mutually singular spectral type.
In view of \cite[Theorem 6.28]{Glasner03} we deduce that they are disjoint in the sense of Furstenberg \cite{Furstenberg67}.

Observe that for any increasing sequence $(N_\ell)$ for which the limit exists, the limit measure
$$\nu:=\lim_{\ell\to\infty}\frac1{N_\ell}\sum_{n=1}^{N_\ell}(R_1\times R_2)^n_*\delta_{(0,y)}$$
on $Y_1\times Y_2$ is a joining of the systems $(Y_1,R_1)$ and $(Y_2,R_2)$.
By disjointness, there only exists the trivial joining and hence
$\nu$ is the product of the two Haar measures $\mu_{Y_1}$ and $\mu_{Y_2}$.
In particular,
\begin{eqnarray*}
  \lim_{N\rightarrow\infty}\frac{1}{N}\sum_{n=1}^NF_1(R_1^n0)F_2(R_2^ny)&=&\int_{Y_1\times Y_2}F_1\otimes F_2\d(\mu_{Y_1}\otimes\mu_{Y_2})\\&=&\int_{Y_1}F_1\d\mu_{Y_1}\cdot\int_{Y_2}F_2\d\mu_{Y_2}=0.
\end{eqnarray*}
This finishes the proof.
\end{proof}

\begin{proof}[Proof of \cref{thm_sharap}]
Let $\theta,\gamma\in\R$ with $\theta>0$ and let $\xbmt$ be a measure preserving system whose discrete spectrum $\sigma(T)$ satisfies $\langle\theta^{-1}\rangle\cap\sigma(T)=\{0\}$. We need to show that
for any $f_1,\ldots,f_k \in L^\infty\xbm$ we have
\begin{equation}\label{eq_proof_thm_sharap}
\lim_{N-M\rightarrow\infty}\frac{1}{N-M}\sum_{n=M}^N \prod_{j=1}^k T^{j\lfloor \theta n+\gamma\rfloor}f_j
=
\lim_{N\rightarrow\infty}\frac{1}{N}\sum_{n=1}^N\prod_{j=1}^k T^{jn}f_j,
\end{equation}
where convergence takes place in $L^2$.

Let us first deal with the case $\theta\geq 1$.
Define $A=\{\lfloor \theta n+\gamma\rfloor:n\in\N\}$
and observe that $m\in A$ if and only if $m\frac1\theta\bmod1\in(a,b]\subset\T$, where $a:=(\gamma-1)/\theta\bmod1$ and $b:=\gamma/\theta\bmod1$.
  Therefore $1_A(m)=1_{(a,b]}(m\frac1\theta)$.
Define
$F(x):=1_{(a,b]}(x)-\frac{1}{\theta}$. Since $F$ is Riemann integrable
with $$\lim_{N\rightarrow\infty}\frac{1}{N}\sum_{m=1}^N F(m\frac{1}
{\theta})=0,$$ it
follows from \cref{lem_hl1} that
$$
\lim_{N-M\rightarrow\infty}\frac{1}{N-M}\sum_{m=M}^N F\left(m\frac1\theta\right)\prod_{j=1}^k T^{jm}f_j
=
0.
$$
From $1_A(m)=F\left(m\frac1\theta\right)+\frac{1}{\theta}$ we deduce that
$$
\lim_{N-M\rightarrow\infty}\frac{1}{N-M}\sum_{m=M}^N 1_A(m) \prod_{j=1}^k T^{jm}f_j
=
\frac{1}{\theta}\left(\lim_{N\rightarrow\infty}\frac{1}{N}\sum_{n=1}^N\prod_{j=1}^k T^{jn}f_j\right).
$$
From this and the observation that the density $d(A)=1/\theta$, \eqref{eq_proof_thm_sharap} follows at once.

Now assume that $\theta<1$.
If $\theta$ is rational, \eqref{eq_proof_thm_sharap} follows immediately, so we assume also that $\theta$ is irrational.
Let $k:=\min\{j\in\N: j\theta>1\}$.
First, we observe that
$\lfloor m\theta+\gamma\rfloor=n$ if and only if
$m\in\left[\lceil{\frac{n-\gamma}{\theta}}\rceil,\lceil{\frac{n+1-\gamma}{\theta}}\rceil\right)$.
Also, the number of $m$'s for which $\lfloor m\theta+\gamma\rfloor=n$ varies between $k-1$ and $k$. Define the sets
$$
A_{k-1}:=\bigg\{n: \Big|\big\{m: \lfloor m\theta+\gamma\rfloor=n\big\}\Big|=k-1\bigg\}=\left\{n:n\frac{1}{\theta}\bmod 1 \in \left(\frac{\gamma}{\theta},1-\xi+\frac{\gamma}{\theta}\right]\right\}
$$
and
$$
A_k:=\bigg\{n: \Big|\big\{m: \lfloor m\theta+\gamma\rfloor=n\big\}\Big|=k\bigg\}=\left\{n:n\frac{1}{\theta}\bmod 1 \in \left(\frac{\gamma}{\theta}-\xi,\frac{\gamma}{\theta}\right]\right\}
$$
and the Riemann integrable functions
$$
F(x):=1_{
\left(\frac{\gamma}{\theta}-\xi,\frac{\gamma}{\theta}\right]}(x)-d(A_k)\qquad\text{and}\qquad
G(x):=1_{\left(\frac{\gamma}{\theta},1-\xi+\frac{\gamma}{\theta}\right]}(x)-d(A_{k-1}).
$$
Observe that $\int F(x)\d x=\int G(x)\d x=0$.
We have
\begin{equation}
\label{eq_proof_thm_sharap-n1}
\begin{split}
\frac{1}{N-M}\sum_{n=M}^N \prod_{j=1}^k T^{j\lfloor \theta n+\gamma\rfloor}f_j=&
\frac{k}{N-M}\sum_{n=\lfloor \theta M+\gamma\rfloor}^{\lfloor \theta N+\gamma\rfloor}1_{A_k}(n)
\prod_{j=1}^k T^{jn}f_j
\\
&+\frac{k-1}{N-M}\sum_{n=\lfloor \theta M+\gamma\rfloor}^{\lfloor \theta N+\gamma\rfloor}1_{A_{k-1}}(n) \prod_{j=1}^k T^{jn}f_j.
\end{split}
\end{equation}
Finally, using $1_{A_{k}}(n)=F(n\frac1\theta)+d(A_{k})$
and $1_{A_{k-1}}(n)=G(n\frac1\theta)+d(A_{k-1})$
we conclude from \cref{lem_hl1} that
$$\frac{1}{N-M}\sum_{n=M}^N \prod_{j=1}^k T^{j\lfloor \theta n+\gamma\rfloor}f_j=\frac{1}{N-M}\sum_{n=M}^N \prod_{j=1}^k T^{jn}f_j.$$
\end{proof}

% =========================================================SECTION
\subsection{Proof of \cref{thm_besicovitch}}
% ================================================================

\begin{proof}[Proof of \cref{thm_besicovitch}]
Let $\phi$ be a Besicovitch almost periodic sequence, let $\xbmt$ be a measure preserving system whose discrete spectrum satisfies $\sigma(T)\cap\langle\sigma(\phi)\rangle=\{0\}$ and let $f_1,\ldots,f_k \in L^\infty(X)$.

We need to show that
$$\lim_{N\rightarrow\infty}\frac{1}{N}\sum_{n=1}^N\phi(n)\prod_{j=1}^k T^{jn}f_j
=
\left(\lim_{N\to\infty}\frac1N\sum_{n=1}^N\phi(n)\right)\left(\lim_{N\rightarrow\infty}\frac{1}{N}\sum_{n=1}^N\prod_{j=1}^k T^{jn}f_j\right).$$

We know from \cite[Theorem 2.24]{Host_Kra09} that all the limits involved exist.
By subtracting from $\phi$ its average,
we can assume without loss of generality that
$ \lim_{N\to\infty}\frac1N\sum_{n=1}^N\phi(n)=0$.
It thus suffices to show that for all $f_0\in L^\infty(X)$ the limit
\begin{equation}\label{eq_proofbesicovitchmulticorrelation}
\lim_{N\rightarrow\infty}\frac{1}{N}\sum_{n=1}^N\phi(n)\int~f_0 \cdot T^{n}f_1\cdot\ldots\cdot  T^{kn}f_k\d\mu=0.
\end{equation}

Using ergodic decomposition if necessary, it suffices to
prove \eqref{eq_proofbesicovitchmulticorrelation}
for ergodic $\mu$.
Invoking \cref{lemma_besicovitchspectrum} we can find for every $\epsilon>0$ a trigonometric polynomial $\rho:n\mapsto\sum_{i=1}^tc_ie(n\theta_i)$ with $\theta_i\in\sigma(\phi)$ that approximates $\phi$ in the sense of \eqref{eq_besicovitchalmostperiodic}.
Observe that $\rho$ can be written as $\rho(n)=F_1(R_1^n0)$
where $R_1:\T^t\to\T^t$ is defined by $R_1(x_1,\dots,x_t)=(x_1+\theta_1,\dots,x_t+\theta_t)$ and $F_1(x_1,\dots,x_t):=\sum_{i=1}^tc_ie(x_i)$ is a continuous function with
$$
\lim_{N\rightarrow\infty}\frac{1}{N}\sum_{n=1}^NF_1(R_1^n0)=0.
$$

In view of \cref{thm:MAIN}, the multi-correlation sequence can be decomposed as
$$\int~f_0\cdot  T^{n}f_1\cdot\ldots\cdot  T^{kn}f_k\d\mu=\psi(n)+\omega(n)+\gamma(n),$$
where $\|\gamma\|_\infty<\epsilon$, $\omega$ is a null-sequence and $\psi(n)=F_2(R_2^ny)$ where $y\in Y$, $F_2\in C(Y)$ and $(Y,R_2)$ is a nilsystem whose discrete spectrum satisfies $\sigma(Y,R_2)\subset\sigma(X,T)$.
Putting everything together we can now rewrite \eqref{eq_proofbesicovitchmulticorrelation} (up to a loss of $2\epsilon$) as
\begin{equation*}
  \lim_{N\rightarrow\infty}\frac{1}{N}\sum_{n=1}^NF_1(R_1^n0)F_2(R_2^ny)=0
\end{equation*}
By arguing as in the proof of \cref{lem_hl1} we see that
the systems $(\T^t,R_1)$ and $(Y,R_2)$ are spectrally disjoint and hence
\begin{eqnarray*}
  \lim_{N\rightarrow\infty}\frac{1}{N}\sum_{n=1}^NF_1(R_1^n0)F_2(R_2^ny)
%&=&\int_{\T^t\times Y}F_1\otimes F_2\d(\mu_{\T^t}\otimes\mu_Y)\\
&=&
\left(\lim_{N\rightarrow\infty}\frac{1}{N}\sum_{n=1}^NF_1(R_1^n0)\right)\left(  \int_YF_2\d\mu_Y\right)
\\
&=&0.
\end{eqnarray*}
This finishes the proof.
\end{proof}

% =========================================================SECTION
\section{Multiple ergodic theorems along primes}\label{sec_primes}
% ================================================================

The purpose of this subsection is to prove
Theorems \ref{thm:along-primes} and \ref{thm_beattyprimes}.
We have a proof which works for both theorems simultaneously and it is presented at the end of this section.
%Note that \cref{thm_beattyprimes} includes \cref{thm:along-primes} as a special case and hence we do not include a separate proof for \cref{thm:along-primes} in this paper.
Our argument for this proof is composed of two key ingredients: The first is
\cref{thm_sharap}, which allows us to control multiple ergodic averages along Beatty sequences. The second ingredient is a variant of a result by Green and Tao regarding the asymptotic Gowers uniformity of the von Mangoldt function. We state this variant below (see \cref{cor:W-trick-beatty}) and include a proof in the appendix.

\subsection{Uniformity norms}

We begin by recalling the definition of the Gowers Uniformity Norms $\|.\|_{U_{[N]}^s}$.
For $N\in\N$ we denote by $[N]$ the set $\{1,2,\ldots,N\}$.

\begin{Definition}[Gowers' Uniformity Norms, \cite{Gowers01}]\label{def:GN}
For $h,N\in\N$ and for $F:\N\rightarrow \C$ let $S^hF$ denote the sequence
$n\mapsto F(n+h)$ and let $F_N$ denote the sequence $F_N(n)=F(N)$ for $n\leq N$ and $F_N(n)=0$ otherwise.
For $s\in\N$ the \define{Gowers Uniformity Norms} $\|.\|_{U^s_{[N]}}$
are defined inductively as
\[
\|F\|_{U^1_{[N]}}:= \left|\frac{1}{N}\sum_{n\in[N]} F_N(n)\right|
\]
and
\[
\|F\|_{U^{s+1}_{[N]}}^{2^{s+1}}:= \frac{1}{N}\sum_{h\in[N]} \left\|F_N S^h \overline{F_N}\right\|_{U^s_{[N]}}^{2^s}.
\]
\end{Definition}

There are different ways of introducing the $U^s_{[N]}$-norms,
which all yield equivalent semi-norms.
%In particular the Gowers uniformity norms $\|\cdot\|_{U^s_{[N]}$ defined above are equivalent to
For comprehensive discussions on
that matter see subsections A.1 and A.2 of Appendix A
in \cite{Frantzikinakis_Host14} or see Appendix B in \cite{Green_Tao10}.

We will make use of the following lemma of Frantzikinakis, Host and Kra.

\begin{Lemma}[Lemma 3.5, \cite{FHK13}]
\label{lem:FHK-Lemma-3.5}
Let $k\in\N$ and let $\xbmt$ be a measure preserving system.
There exists $M\in\N$ such that for any sequence $F:\N\rightarrow \C$
satisfying $\frac{F(n)}{n^c}\to 0$ for all $c>0$, and any $f_0,f_1,\ldots,f_k \in L^\infty\xbm$ we have
\begin{equation}
\label{eqn:FHK-Lemma-3.5-1}
\left\|\frac{1}{N}\sum_{n=1}^N F(n) \prod_{j=1}^k T^{jn}f_j\right\|_{L^2\xbm}
\leq  M\|F\|_{U^k_{[N]}} +o_N(1).
\end{equation}
\end{Lemma}

\subsection{A variant of the Green-Tao-Ziegler theorem}
Let $\Lambda:\N\to\R$ denote the von Mangoldt function, defined as $\Lambda(p^k)=\log(p)$ for every prime $p$ and $k\geq0$, and $\Lambda(n)=0$ otherwise. For $b,W\in\N$
define
$$
\Lambda'(n):=\log(n)1_\P(n)\qquad\text{and}\qquad
\Lambda_{W,b}(n):=\frac{\phi(W)}{W}\Lambda'(Wn+b),
$$
where $\phi$ denotes Euler's totient function.

In \cite{Green_Tao10} Green and Tao prove the following result (initially conditional on two conjectures, which were eventually settled by them and Ziegler \cite{GT12-2,GTZ12}).

\begin{Theorem}[\cite{Green_Tao10,GT12-2,GTZ12}; also, cf. Theorem 2.2 in \cite{FHK13}]
\label{thm:W-trick}
Let $s,w>1$ and let $W:=\prod_{p\leq w}p$ be the product of all primes up to $w$. Then
\[
\lim_{w\to\infty}\limsup_{N\to\infty}\max_{\substack{b\leq W\\(b,W)=1}}\|\Lambda_{W,b}-1\|_{U^s_{[N]}}=0
\]
%converges to $0$ as $N\to\infty$ and then $w\to\infty$.
\end{Theorem}

Next, consider the following modified versions of the von Mangoldt function. Let $\theta,\gamma\in\R$ with $\theta>0$ and define
$$
\Lambda_{\theta,\gamma}(n):=\log(n) 1_{\P(\theta,\gamma)}(n)\qquad\text{and}\qquad
\Lambda_{\theta,\gamma,W,b}(n):=\frac{\phi(W)}{W}\Lambda_{\theta,\gamma}(Wn+b),
$$

The following extension of \cref{thm:W-trick} follows from part (2) of Proposition 3.2 in \cite{Sun15}.
\begin{Theorem}
\label{cor:W-trick-beatty}
Let $s,w>1$, let $\theta,\gamma\in\R$ with $\theta>0$ irrational and let $W:=\prod_{p\leq w}p$ be the product of all primes up to $w$.
Let $\B(\theta,\gamma,W,b):=\{n\in\N: Wn+b\in \{\lfloor\theta m+\gamma\rfloor:m\in\Z\}\}$.
Then
\[
\lim_{w\to\infty}\limsup_{N\to\infty}\max_{\substack{b\leq W\\(b,W)=1}}\|\Lambda_{\theta,\gamma,W,b}-1_{\B(\theta,\gamma,W,b)}\|_{U^s_{[N]}}=0.
\]
\end{Theorem}

\subsection{Proof of \cref{thm_beattyprimes}}

We employ the following folklore lemma to transfer information about the von Mangoldt function to information about the primes (see for instance \cite[Lemma 1]{Frantzikinakis_Host_Kra07}).
\begin{Lemma}\label{lemma_Hilbertprime}Let $H$ be a Hilbert space and let $\alpha:\Z\to H$ be a bounded sequence.
Then
  $$\lim_{N\rightarrow\infty}\frac{1}{\pi(N)}\sum_{p\leq N }
\alpha(p)=
\lim_{N\rightarrow\infty}\frac{1}{N}\sum_{n=1}^N \Lambda'(n)\alpha(n)$$
in the sense that if one of the limits exists (in the strong topology) then so does the other and they are equal.
\end{Lemma}

From \cref{lemma_Hilbertprime} we immediately obtain the following corollary regarding
Beatty sequencers in primes.

\begin{Corollary}
\label{lemma_Hilbertbeattyprime}
Let $H$ be a Hilbert space and let $\alpha:\Z\to H$ be a bounded sequence
and let $\pi_{\theta,\gamma}(N)$ denote the number of primes in $\P(\theta,\gamma)$ up to $N$.
Then
  $$\lim_{N\rightarrow\infty}\frac{1}{\pi_{\theta,\gamma}(N)}\sum_{p\in [N]\cap \P(\theta,\gamma)}
\alpha(p)=
\lim_{N\rightarrow\infty}\frac{\theta}{N}\sum_{n=1}^N \Lambda_{\theta,\gamma}(n)\alpha(n)$$
in the sense that if one of the limits exists (in the strong topology) then so does the other and they are equal.
\end{Corollary}

\cref{lemma_Hilbertbeattyprime} follows from applying \cref{lemma_Hilbertprime} to $\tilde\alpha(n)=\alpha(n)1_{\{\lfloor\theta m+\gamma\rfloor:m\in\Z\}}(n)$, together with the Prime Number Theorem for Beatty sequences (cf. \cite[end of page 289]{Ribenboim96}), which states that
$\pi_{\theta,\gamma}(N)\sim \frac{N}{\theta \log(N)}$.

We are now ready to give a proof of \cref{thm_beattyprimes}.

\begin{proof}[Proof of \cref{thm_beattyprimes}]
Let $\theta,\gamma\in\R$ be such that either $\theta=1$ or $\theta$ is positive and irrational, let $\P(\theta,\gamma)=\P\cap\{\lfloor n\theta+\gamma\rfloor:n\in\Z\}$ and let $\pi_{\theta,\gamma}(N)=\big|\{1,\dots,N\}\cap\P(\theta,\gamma)\big|$.
Let $\xbmt$ be a measure preserving system whose discrete spectrum $\sigma(T)$ satisfies $\sigma(T)\cap\langle\Q,\theta^{-1}\rangle=\{0\}$, let $f_1,\ldots,f_k \in L^\infty\xbm$ and define
$$\alpha(n):=\prod_{j=1}^kT^{jn}f_j.$$
We want to show that
$$\lim_{N\to\infty}\frac1{\pi_{\theta,\gamma}(N)}\sum_{\substack{p\leq N,\\p\in\P(\theta,\gamma)}}~
\alpha(p)
=
\lim_{N\rightarrow\infty}\frac{1}{N}\sum_{n=1}^N~\alpha(n).$$

We remark that the case $\theta\in(0,1)$ follows trivially from the case $\theta=1$, so we can assume without loss of generality that $\theta\geq1$.

In view of \cref{lemma_Hilbertbeattyprime} it remains to show that
$$
\lim_{N\to\infty}\frac{\theta}{N}\sum_{n=1}^N\Lambda_{\theta,\gamma}(n)\alpha(n)=
\lim_{N\to\infty}\frac1N\sum_{n=1}^N\alpha(n).
$$

Let $\epsilon>0$ be arbitrary.
Using \cref{cor:W-trick-beatty}, for large enough $w,N\in\N$, letting $W:=\prod_{p\leq w}p$, we have
\begin{equation}\label{eq_Wtrick}
\max_{\substack{b\leq W,\\(b,W)=1}}\|\Lambda_{\theta,\gamma,W,b}(n)-1_{\B(\theta,\gamma,W,b)}\|_{U^s_{[N]}}\leq \epsilon.
\end{equation}
For convenience, assume that $N$ is a multiple of $W$.
We have:
\begin{eqnarray*}
\frac{\theta}{N}\sum_{n=1}^N \Lambda_{\theta,\gamma}(n)\alpha(n)
&=&\frac{\theta}{N}\sum_{n=1}^N
\sum_{\substack{b\leq W,\\(b,W)=1}} 1_{W\Z+b}(n)\Lambda_{\theta,\gamma}(n)\alpha(n)+o_N(1)\\
&=&\sum_{\substack{b\leq W,\\(b,W)=1}}\frac{\theta}{N}\sum_{n=0}^{N/W}\Lambda_{\theta,\gamma}(Wn+b)\alpha(Wn+b)+o_N(1)\\
&=&\frac{W}{\phi(W)}\sum_{\substack{b\leq W,\\(b,W)=1}}\frac {\theta}{N}\sum_{n=0}^{N/W}\Lambda_{\theta,\gamma,W,b}(n)\alpha(Wn+b)+o_N(1).
\end{eqnarray*}
Next, we fix $b\leq W$ satisfying $(b,W)=1$ and decompose $\Lambda_{\theta,\gamma,W,b}(n)=
\big(\Lambda_{\theta,\gamma,W,b}(n)-1_{\B(\theta,\gamma,W,b)}\big)+1_{\B(\theta,\gamma,W,b)}$.
%Observe that, for each $i$, the map $n\mapsto p_i(Wn+b)$ is a polynomial of the same degree as $p_i$.

According to \cref{lem:FHK-Lemma-3.5}, equation \eqref{eqn:FHK-Lemma-3.5-1} holds for all $F:\N\to\C$
that satisfy $\frac{F(n)}{n^c}\to 0$ for all $c>0$; it is clear that $F(n):=\Lambda_{\theta,\gamma,W,b}(n)-1_{\B(\theta,\gamma,W,b)}$
satisfies this condition for every $W,b\in\N$, where $\B(\theta,\gamma,W,b):=\{n\in\N: Wn+b\in \{\lfloor\theta m+\gamma\rfloor:m\in\Z\}\}$.
Therefore, in view of \eqref{eqn:FHK-Lemma-3.5-1} and \eqref{eq_Wtrick}, we get
\begin{equation*}
\begin{split}
\frac{\theta}{N}\sum_{n=0}^{N/W}\Lambda_{\theta,\gamma,W,b}(n)&\alpha(Wn+b)
\\
&=\frac{\theta }N\sum_{n=0}^{N/W}1_{\B(\theta,\gamma,W,b)}(n)\alpha(W n + b)+R(N)
\\
&=\frac{\theta}{N}\sum_{n=0}^{N} 1_{\{\lfloor \theta m+\gamma\rfloor:m\in\Z\}}(n)1_{W\Z+b}(n) \alpha(n)+R(N)
\end{split}
\end{equation*}
where $\|R(N)\|_{L^2}\leq M\|\Lambda_{\theta,\gamma,W,b}-1_{\B(\theta,\gamma,W,b)}\|_{U^s_{[N]}}+o_N(1)\leq M\epsilon+o_N(1)$ for some $M\in\R$ which does not depend on $\epsilon$.
%\begin{eqnarray*}
%\frac{1}{N}\sum_{n=1}^N\Lambda_{Wn+b}'(n)\varphi(Wn+b)\\
%&=&\frac{1}{\phi(W)}\sum_{\substack{b\leq W,\\(b,W)=1}}
%\frac{1}{N}\sum_{n=1}^N\varphi(Wn+b).
%\end{eqnarray*}

According to \cref{lemma_productbesicovitch}, the product $1_{\{\lfloor \theta m+\gamma\rfloor:m\in\Z\}}\cdot1_{W\Z+b}$ is Besicovitch almost periodic whose spectrum is contained in $\langle \Q\cup\{\theta^{-1}\}\rangle$.
Since the Ces\`aro average of the product of two Besicovitch almost periodic sequences with disjoint spectra is the product of the Ces\`aro averages of the two sequences, it follows that
$$\lim_{N\to\infty}\frac1N\sum_{n=1}^N1_{\{\lfloor \theta m+\gamma\rfloor:m\in\Z\}}(n)\cdot1_{W\Z+b}(n)=\frac1{\theta W}$$
for every $\theta$ which is either irrational or co-prime with $W$.

Using \cref{thm_besicovitch} we thus obtain
\[
\lim_{N\rightarrow\infty}\frac{\theta}{N}\sum_{n=1}^N 1_{\{\lfloor \theta m+\gamma\rfloor:m\in\Z\}}(n)1_{W\Z+b}(n)
\alpha(n) =
\lim_{N\rightarrow\infty}\frac{1}{WN}\sum_{n=1}^{N}\alpha(n)
\]
for all $W$ and $b$. Therefore, averaging over all $b$ with $(b,W)=1$ we get
\[
\limsup_{N\rightarrow\infty}\left\|\frac{\theta}{N}\sum_{n=1}^N \Lambda_{\theta,\gamma}(n) \alpha(n)-\frac{1}{N}\sum_{n=1}^N\alpha(n)\right\|_{L^2}\leq M\epsilon.
\]
Since $\epsilon>0$ was arbitrary, this finishes the proof.
\end{proof}

% =========================================================SECTION
\section{Spectrum of the orbit of the diagonal}
\label{sec_main-tehnical-result}
% ================================================================

The purpose of this section is to give a proof of
\cref{cor:thm3.4}.
For convenience, let us restate the theorem here.

\begin{named}{\cref{cor:thm3.4}}{}
Let $k\in\N$, let $X$ be a connected nilmanifold and let $R:X\to X$
be an ergodic nilrotation.
Define $S:= R\times R^2\times \ldots\times R^k$ and
\begin{equation}\label{eq_cor:thm3.4}
Y_{X^\Delta}:=\overline{\big\{S^n(x,x,\ldots, x):
x\in X,~n\in\Z\big\}}\subset X^{k}.
\end{equation}
Then $\sigma(X,R)=\sigma(Y_{X^\Delta}, S)$.
\end{named}

We will derive \cref{cor:thm3.4} from the following more general
result.

\begin{Theorem}\label{thm:3.4second}
Let $k\in\N$, let $X$ be a connected nilmanifold and let $R:X\to X$
be an ergodic nilrotation.
Define $S:= R\times R^2\times \ldots\times R^k$ and
\begin{equation}
\label{eq:yx}
Y_x:=\overline{\big\{S^n(x,x,\ldots, x): n\in\Z\big\}}
\subset X^{k}.
\end{equation}
Then for almost every $x\in X$ the Kronecker factor
of $(Y_x,S)$ is isomorphic to the Kronecker factor $(K,R)$ of $(X,R)$.
In particular, $\sigma(X,R)=\sigma(Y_x, S)$ for almost every $x\in X$.
\end{Theorem}

We remark that in the statement of
\cref{thm:3.4second} the restriction of $x$ to a full measure subset of $X$ is necessary, because
there may be points $x\in X$ for which the spectrum of the system
$(Y_x, S)$ is strictly larger than the spectrum of $(X,R)$, as the following example illustrates.

\begin{Example}
Consider the matrix group
$$
G:=\left\{
\begin{pmatrix}
 1&n&c\\
 0&1&b\\
 0&0&1
\end{pmatrix}
: n\in\Z,~b,c\in\R\right\}.
$$
This group is a $2$-step nilpotent Lie group and it acts continuously and transitively on $\T^2$
via
$$
\begin{pmatrix}
 1&n&c\\
 0&1&b\\
 0&0&1
\end{pmatrix}(y,z) = (y+b,z+c+ny),\qquad \forall (y,z)\in\T^2.
$$
Also, $(\T^2,G)$ is a nilsystem since it is isomorphic to
$(G/\Gamma,G)$, where $\Gamma$ is the uniform and discrete subgroup
of $G$ given by
$$
\Gamma:=\left\{
\begin{pmatrix}
 1&n&m\\
 0&1&k\\
 0&0&1
\end{pmatrix}
: n,m,k\in\Z\right\}.
$$
Let $\alpha$ be an arbitrary irrational number and let
$a\in G$ denote the element
$$
a=
\begin{pmatrix}
 1&1&0\\
 0&1&\alpha\\
 0&0&1
\end{pmatrix}.
$$
Then the nilrotation $R_a:\T^2\to\T^2$, which takes the from $(y,z)\mapsto (y+\alpha, z+y)$, is totally ergodic.
However, if $x$ denotes the point $x:=(\tfrac14,0)$, then the closure $Y$ of the set
$$
\{(R_a^nx,R_a^{2n}x):n\in\Z\}
$$
in $\T^2\times\T^2$ is not connected. In fact, straightforward
calculations reveal that it consists of two connected components.
This implies that $1/2\in\sigma(Y,R_a\times R_{a^2})$ but $1/2\notin\sigma(\T^2,R_a)$.
\end{Example}

The proof of \cref{thm:3.4second} is presented in Sections \ref{sec_connected-x} through
\ref{sec_spec-case-connected}.
For now we will present the deduction of \cref{cor:thm3.4} from \cref{thm:3.4second}.

\begin{proof}[Proof of \cref{cor:thm3.4} assuming \cref{thm:3.4second}]
Let $k\in\N$, let $X$ be a connected nilmanifold and let $R:X\to X$
be an ergodic nilrotation.
Define $S:= R\times R^2\times \ldots\times R^k$ and let $Y_{X^\Delta}$ be defined by \eqref{eq_cor:thm3.4}.
Given $\theta\in\sigma(X,R)$, let $f\in L^2(X)$ be an eigenfunction of the system $(X,R)$ with eigenvalue $\theta$.
Since the function $\tilde f\in L^2(Y_{X^\Delta})$ defined by $\tilde f(x_1,\dots,x_k)=f(x_1)$ is an eigenfunction for the system $(Y_{X^\Delta},S)$ with eigenvalue $\theta$, it follows that $\sigma(X,R)\subset\sigma(Y_{X^\Delta},S)$.

Next we prove the converse inclusion.
Let $\nu$ be the Haar measure of the nilmanifold $Y_{X^\Delta}$ and let $\nu_x$ be the Haar measure of the nilmanifold $Y_x$ defined by \eqref{eq:yx}.
Observe that the sets $Y_x$ are precisely the atoms of the invariant $\sigma$-algebra of the system $(Y_{X^\Delta},S)$.
Therefore, the measures $\nu_x$ form the ergodic decomposition of $\nu$.

Let $\theta\in\sigma(Y_{X^\Delta},S)$ and let $f\in L^2(Y_{X^\Delta},\nu)$ be an eigenfunction with eigenvalue $\theta$.
In other words, for almost every $y\in Y_{X^\Delta}$ we have $Sf(y)=e(\theta)f(y)$.
Since $f$ cannot be $0$ $\nu$-a.e., there exists a positive measure set of $x\in X$ for which the restriction of $f$ to the system $(Y_x,\nu_x,S)$ is not the zero function.
But for any such $x$, the restriction of $f$ to the system $(Y_x,\nu_x,S)$ is an eigenfunction with eigenvalue $\theta$.
Finally, we invoke \cref{thm:3.4second} to conclude that $\theta\in\sigma(X,R)$, finishing the proof.
\end{proof}

% =========================================================SECTION
\subsection{A useful reduction}%At this point $X$ is always connected... Perhaps $G$?
\label{sec_connected-x}
% ================================================================

Let $G$ be an $s$-step nilpotent Lie group, let $\Gamma$
be a uniform and discrete subgroup of $G$ and assume that
$X=G/\Gamma$ is connected. This is equivalent to the assertion
that $G=G^\circ \Gamma$, where $G^\circ$ is the connected
component of the identity in $G$.

Let $a\in G$ be arbitrary and consider the group
$G':=\langle G^\circ, a\rangle$ generated by $G^\circ$ and $a$.
We have the following useful lemma regarding $G'$.

\begin{Lemma}
\label{lem:group-generated-by-identity-component-end-element}
The group $G'$ is a closed rational subgroup of $G$ and hence a $s$-step
nilpotent Lie group.
If
\[
G'=G'_1 \trianglerighteq G'_2\trianglerighteq G'_3\trianglerighteq \ldots\trianglerighteq G'_{s}\trianglerighteq\{1_G'\}
\]
denotes the lower central series of $G'$, then
$G'_i$ is connected for $2\leq i\leq s$.
\end{Lemma}

\begin{proof}
First, let us show that $G'$ is closed.
Since $G^\circ$ is a clopen subset of $G$ we deduce that $G'$ is a disjoint union of clopen sets and therefore clopen. In particular, $G'$ is closed.

Next, we show that $G'_i$ is connected for $2\leq i\leq s$.
Since $G^\circ$ is a normal subgroup of $G$ it holds
that $\langle G^\circ, a\rangle= G^\circ a^\Z$. Therefore
$$
G'_2=[G^\circ a^\Z,G^\circ a^\Z]=
[G^\circ,G^\circ ][G^\circ, a^\Z][a^\Z,a^\Z]=
[G^\circ,G^\circ ][G^\circ ,a^\Z].
$$
Note that groups generated by connected sets are connected.
Hence the group $[G^\circ,G^\circ ][G^\circ ,a^\Z]$ is connected, as it is generated by the connected set
$$
\bigcup_{g\in G^\circ}[G^\circ,g]\cup\bigcup_{n\in\Z}[G^\circ,a^n].
$$
Hence, $G'_2$ is connected.
Analogous arguments can be used to show by induction on $i$
that  $G'_i$ is connected for all $i=3,\ldots, s$.

Finally, to see why $G'$ is a rational subgroup of $G$% whenever
%$G/\Gamma$ is connected
, let $\pi:G\to G/\Gamma$ denote the natural projection from $G$ onto $G/\Gamma$.
In view of \cref{rem:rational-subgroup-equivalences}, $G'$ is rational
if and only if $\pi(G')$ is closed. Since $G/\Gamma$ is connected
it follows that $\pi(G^\circ)=G/\Gamma$ and so it follows from
the fact that $G^\circ$ is contained in $G'$ that
$\pi(G')=G/\Gamma$ is closed. This finishes the proof.
\end{proof}

\begin{Remark}\label{rem_gen-conn-el}
Let $(X,G, R_a)$ be an ergodic nilsystem
with connected phase space $X=G/\Gamma$
and let $G':=\langle G^\circ,a\rangle$ be the group generated by
$G^\circ$ and $a$.
By \cref{lem:group-generated-by-identity-component-end-element}
the group $G'$ is rational which means that
$\Gamma':=\Gamma\cap G'$ is a uniform and discrete subgroup of
$G'$. Let $X'$ denote the nilmanifold $G'/\Gamma'$ and consider the
nilsystem $(X',G', R_a)$.

We claim that the two nilsystems
$(X,G, R_a)$ and $(X',G', R_a)$ are isomorphic. To verify this
claim consider the map $\eta:X'\to X$ defined by the formula
$\eta(g\Gamma')=g\Gamma$ for all $g\in G'$.
This map is well defined, continuous and injective.
Moreover, using the fact
that $G=G^\circ \Gamma$ and that $G^\circ\subset G'$, we conclude
that $\eta$ is also surjective. Since any continuous bijection
between compact Hausdorff spaces is a homeomorphism we get
that $X'$ and $X$ are indeed homeomorphic spaces. Finally, since
$R_a\circ \eta=\eta\circ R_a$, we conclude that
$(X,G, R_a)$ and $(X',G', R_a)$ are isomorphic dynamical systems.
\end{Remark}
% =========================================================SECTION
\subsection{The subgroup $H$}
\label{sec_subgroup-H}
% ================================================================

Let $G$ be an $s$-step nilpotent Lie group, let $\Gamma$ be a uniform and
discrete subgroup such that $X=G/\Gamma$ is connected and let
$R=R_a$ be an ergodic nilrotation.
It will be convenient for us to assume that $G$
is generated by $G^\circ$ and $a$. Note that in view of
\cref{rem_gen-conn-el} this assumption can be made without loss
of generality.
Let $G=G_1 \trianglerighteq G_2\trianglerighteq
\ldots\trianglerighteq G_{s}\trianglerighteq\{1_G\}$
denote the lower central series of $G$ and
fix some positive integer $k\in\N$. We define the subsets
$H^{(1)},\ldots,H^{(k-1)}$ of $G^k$ as
\begin{equation}\label{eqn:h-i}
H^{(i)}:=\{
(g^{\binom{1}{i}},
g^{\binom{2}{i}},\ldots,g^{\binom{k}{i}}): g\in G_i\},
\end{equation}
where $\binom{j}{i}=0$ for $j<i$, and we define $H$ as
\begin{equation}\label{eqn:h}
H:= H^{(1)} H^{(2)}\cdots H^{(k-1)} G_k^k.
\end{equation}
The set $H$ is in fact a subgroup of $G^{k}$
and can be used to explicitly describe
the orbit closure of the diagonal in $X^{k}$ under the
nilrotation $S:=R_u\times R_u^2\times \ldots\times R_u^k$.
The next proposition lists some of the well known properties of $H$;
for a more comprehensive discussion
on $H$ see \cite{Ziegler05,BHK05,Leibman10b}.

\begin{Proposition}
\label{prop:H}
Let $G$ be a $s$-step nilpotent Lie group,
let $\Gamma$ be a uniform and discrete subgroup of $G$
such that $X=G/\Gamma$ is connected, let $a\in G$ and
assume that $G=\langle G^\circ, a\rangle$. Let $k\in\N$ and
define $H$ as in \eqref{eqn:h}. Then
\begin{enumerate}	
[label=(\roman*),ref=(\roman*),leftmargin=*]
\item\label{prop:H-1}
$H$ is a closed and rational subgroup of $G^k$.
%generated by $H^\circ$
%and $(1_G,u,u^2,\ldots,u^k)$.
\item\label{prop:H-2}
The commutator subgroup $[H,H]$ of $H$ satisfies
$$
[H,H]= H\cap [G^k,G^k].% = [H_1,H_1]H_2\cdots H_{k-1}G_k^k.
$$
\item\label{prop:H-3}
Define $\Delta:=H\cap \Gamma^k$. Then $\Delta$
is a uniform and discrete subgroup of $H$ and
$Y:=H/\Delta$ is a connected nilmanifold.
\item\label{prop:H-4}
For $b\in G$ define
$$
S_b:= R_b\times R_b^2\times \ldots\times R_b^k
$$
and
for $x=g\Gamma \in X$ define
$$
Y_x:=\overline{\big\{S_a^n(x,x,\ldots, x): n\in\Z\big\}}\subset X^k.
$$
For almost every $x=g\Gamma\in X$
the nilsystems $(Y_x, S_a)$ and $(Y,S_{g^{-1}ag})$ are isomorphic.
\end{enumerate}
\end{Proposition}

\begin{proof}
For the proofs of items \ref{prop:H-1} -- \ref{prop:H-3} we refer the reader to
\cite{Leibman98}, \cite[Theorem 5.1]{BHK05} and \cite[Proposition 5.7]{Leibman10b}.

For the proof of \ref{prop:H-4} we repeat a short argument
that appeared in \cite[Subsection 2.5]{Frantzikinakis08}.
We need the following theorem.

\begin{Theorem}[see {\cite[Theorem 2.2]{Ziegler05}} or {\cite[Theorem 5.4]{BHK05}}]
Let $G$ be a $s$-step nilpotent Lie group,
let $\Gamma$ be a uniform and discrete subgroup of $G$
such that $X=G/\Gamma$ is connected, let $a\in G$ be such that $R_a$
is ergodic and assume that $G=\langle G^\circ, a\rangle$. Let $k\in\N$,
define $H$ as in \eqref{eqn:h} and put
$\Delta:=H\cap \Gamma^k$ and
$Y:=H/\Delta$. If $f_1,\ldots,f_k \in L^\infty (X)$, then for
a.e. $x=g\Gamma \in X$ we have
\begin{equation}
\label{eq:zlf}
\lim_{N\to\infty}\frac{1}{N}\sum_{n=1}^N
f_1(g^{-1}a^n g\Gamma)\cdot\ldots\cdot f_k({g^{-1}a^{kn}g}\Gamma)
=
\int_{Y} f_1\otimes\ldots\otimes f_k~\d\mu_{Y}.
\end{equation}
\end{Theorem}

In the following let us identify $Y$ with its embedding into $X^k$
in the obvious manner. Consider the injective map
$R_{g^{-1}}\times\ldots\times R_{g^{-1}}:X^k\to X^k$.
It follows from the definition of the
group $H$ that for any $n\in\Z$ the image of the point
$S_a^n(x,\ldots,x)$ under $R_{g^{-1}}\times\ldots\times R_{g^{-1}}$
lies in $Y$. Since points of the form
$S_a^n(x,x,\ldots, x)$ are dense in $Y_x$ we conclude that
$$
(R_{g^{-1}}\times\ldots\times R_{g^{-1}})(Y_x)\subset Y.
$$
Hence, it only remains to show that
for almost all $x\in X$ we have
$(R_{g^{-1}}\times\ldots\times R_{g^{-1}})(Y_x)=Y$. This, however,
follows right away from the fact that
$\{S_{g^{-1}ag}^n(x,x,\ldots, x): n\in\Z\}$ is contained in
$(R_{g^{-1}}\times\ldots\times R_{g^{-1}})(Y_x)$ and that
\eqref{eq:zlf} implies
$\{S_{g^{-1}ag}^n(x,x,\ldots, x): n\in\Z\}$ is dense in $Y$
for almost all $x\in X$.
\end{proof}

% =========================================================SECTION
\subsection{A proof of \cref{thm:3.4second} }
\label{sec_spec-case-connected}
% ================================================================

\begin{proof}[Proof of \cref{thm:3.4second}]
Let $k\in\N$, let $X$ be a connected nilmanifold and let $R=R_a$
be an ergodic nilrotation. In view of
\cref{rem_gen-conn-el} we can assume without loss of
generality that $G$ is generated by $G^\circ$ and $a$.
Let $S=S_a:=R\times R^2\times\ldots\times R^k$,
and let $(K,R)$ denote the Kronecker factor of $(X,R)$.
We want to show that for almost every $x\in X$ the Kronecker factor
of $(Y_x,S)$ is isomorphic (as measure preserving system) to the system $(K,R)$, where $Y_x$ is defined by \eqref{eq:yx}.

Let $H$ be as in \eqref{eqn:h}.
Following the notation of
\cref{prop:H}, we define $Y:=H/\Delta$, where
$\Delta:=H\cap \Gamma^k$, and
$S_{b}:=R_{b}\times R_{b}^2\times \ldots\times R_{b}^k$.
By \cref{prop:H} part \ref{prop:H-4} the system
$(Y_x, S_a)$ is isomorphic to $(Y,S_{g^{-1}a g})$,
for almost every $x=g\Gamma\in X$.
Hence, to finish the proof
it suffices to show that the Kronecker factor of
$(Y,S_{b})$ is isomorphic to $(K,R)$
for every $b$ of the form $g^{-1}ag$.

Let $H^{(1)},\ldots, H^{(k-1)}$ be as in \eqref{eqn:h-i}.
According to
\cref{lem:group-generated-by-identity-component-end-element}
the groups $G_2,G_3,\ldots,G_s$ are connected, from which we deduce
that the group $G_k^k$ and the sets $H^{(2)},\ldots, H^{(k-1)}$ are also
connected.
Since $G$ is generated by $G^\circ$ and $a$, we get that the group generated by
$H^{(1)}$ is generated by its identity component and $(a,a^2,\ldots,a^k)$.
Putting everything together, this implies that
$H$ is generated by $H^\circ$ and $(a,a^2,\ldots,a^k)$.
Moreover, since
$H$ is normalized
by the diagonal $G^\Delta$, we conclude that $H$ is
in fact generated by $H^\circ$ and
$h_g:=(g^{-1}ag,g^{-1}a^2g,\ldots,g^{-1}a^kg)$
for any $g\in G$.

Now, the fact that $H$ is generated by $H^\circ$
and $h_g$ allows us to apply \cref{cor:Thm2.17-kronecker} to the system $(Y,S_b)$, noting that this system is isomorphic to systems of the form $(Y_x,S_a)$ which are, by construction, transitive and hence, in view of \cref{thm:dynamics-nilrotation}, ergodic.
It now follows that for any $b=g^{-1}ag$ the
Kronecker factor of $(Y,S_b)$ is given by
$(K_Y, S_b)$, where $N_Y=[H,H]$ and $K_Y:=N_Y\backslash Y$.

Finally, we need to show that the systems $(K_Y,S_b)$ and $(K,R)$ are isomorphic.
Observe that $K_Y$ can be identified with $H/(N_Y\Delta)$.
Moreover, by definition we have $\Delta = H\cap \Gamma^k$
and from \cref{prop:H}, part \ref{prop:H-2}, we have $N_Y= H\cap N^k$, where $N=[G,G]$.
This implies that
$N_Y\Delta=H\cap (N^k\Gamma^k)$ and hence $K_Y=H/(H\cap N^k\Gamma^k)$.
It follows from the second isomorphism theorem that $K_Y\cong (HN^k\Gamma^k)/(N^k\Gamma^k)$.
Finally, since $H^{(i)}\subset N^k$ for every $i>1$, it follows that $HN^k=H^{(1)}N^k$ and thus
$K_Y\cong (H^{(1)}N^k\Gamma^k)/(N^k\Gamma^k)$.

Applying \cref{cor:Thm2.17-kronecker} to the original system $(X,R)$, we see that $K=N\backslash X=G/N\Gamma$; therefore $K_Y$ embeds naturally into $K^k=G^k/(N^k\Gamma^k)$.
Looking at the definition of $H^{(1)}$ we deduce that
$$
K_Y\cong \{(v,2v,\ldots,kv):v\in K\}\cong K.
$$
Finally, since $K$ is an abelian group and $b=g^{-1}ag$, the projections of $a$ and $b$ onto $K$ coincide. Hence $(K_Y,S_b)$ is indeed isomorphic (as a measure preserving system) to $(K,R_a)$ as desired.
\end{proof}

\appendix

% ====================================================BIBLIOGRAPHY

\bibliographystyle{siam}

\providecommand{\noopsort}[1]{}
\allowdisplaybreaks
\small
\bibliography{BibMR15,refs-joel}

\begin{thebibliography}{10}

\bibitem{AGH63}
{\sc L.~Auslander, L.~Green, and F.~Hahn}, {\em Flows on homogeneous spaces},
  With the assistance of L. Markus and W. Massey, and an appendix by L.
  Greenberg. Annals of Mathematics Studies, No. 53, Princeton University Press,
  Princeton, N.J., 1963.

\bibitem{BL85}
{\sc A.~Bellow and V.~Losert}, {\em The weighted pointwise ergodic theorem and
  the individual ergodic theorem along subsequences}, Trans. Amer. Math. Soc.,
  288 (1985), pp.~307--345.

\bibitem{BHK05}
{\sc V.~Bergelson, B.~Host, and B.~Kra}, {\em Multiple recurrence and
  nilsequences}, Invent. Math., 160 (2005), pp.~261--303.
\newblock With an appendix by Imre Ruzsa.

\bibitem{Besicovitch55}
{\sc A.~S. Besicovitch}, {\em Almost periodic functions}, Dover Publications,
  Inc., New York, 1955.

\bibitem{Bohr25}
{\sc H.~Bohr}, {\em Zur theorie der fast periodischen funktionen}, Acta Math.,
  45 (1925), pp.~29--127.
\newblock I. Eine verallgemeinerung der theorie der fourierreihen.

\bibitem{Bohr25b}
\leavevmode\vrule height 2pt depth -1.6pt width 23pt, {\em Zur {T}heorie der
  {F}astperiodischen {F}unktionen}, Acta Math., 46 (1925), pp.~101--214.
\newblock II. Zusammenhang der fastperiodischen Funktionen mit Funktionen von
  unendlich vielen Variabeln; gleichm{\"a}ssige Approximation durch
  trigonometrische Summen.

\bibitem{Durrett10}
{\sc R.~Durrett}, {\em Probability: theory and examples}, Cambridge Series in
  Statistical and Probabilistic Mathematics, Cambridge University Press,
  Cambridge, fourth~ed., 2010.

\bibitem{Frantzikinakis04}
{\sc N.~Frantzikinakis}, {\em The structure of strongly stationary systems}, J.
  Anal. Math., 93 (2004), pp.~359--388.

\bibitem{Frantzikinakis08}
{\sc N.~Frantzikinakis}, {\em Multiple ergodic averages for three polynomials
  and applications}, Trans. Amer. Math. Soc., 360 (2008), pp.~5435--5475.

\bibitem{Frantzikinakis_Host14}
{\sc N.~Frantzikinakis and B.~Host}, {\em Higher order fourier analysis of
  multiplicative functions and applications}.
\newblock To appear in J. Amer. Math. Soc.; available online at
  http://arxiv.org/abs/1403.0945.

\bibitem{Frantzikinakis_Host_Kra07}
{\sc N.~Frantzikinakis, B.~Host, and B.~Kra}, {\em Multiple recurrence and
  convergence for sequences related to the prime numbers}, J. Reine Angew.
  Math., 611 (2007), pp.~131--144.

\bibitem{FHK13}
{\sc N.~Frantzikinakis, B.~Host, and B.~Kra}, {\em The polynomial
  multidimensional {S}zemer\'edi theorem along shifted primes}, Israel J.
  Math., 194 (2013), pp.~331--348.

\bibitem{Furstenberg67}
{\sc H.~Furstenberg}, {\em Disjointness in ergodic theory, minimal sets, and a
  problem in {D}iophantine approximation}, Math. Systems Theory, 1 (1967),
  pp.~1--49.

\bibitem{Furstenberg77}
\leavevmode\vrule height 2pt depth -1.6pt width 23pt, {\em Ergodic behavior of
  diagonal measures and a theorem of {S}zemer\'edi on arithmetic progressions},
  J. Analyse Math., 31 (1977), pp.~204--256.

\bibitem{Glasner03}
{\sc E.~Glasner}, {\em Ergodic theory via joinings}, vol.~101 of Mathematical
  Surveys and Monographs, American Mathematical Society, Providence, RI, 2003.

\bibitem{Gowers01}
{\sc W.~T. Gowers}, {\em A new proof of {S}zemer\'edi's theorem}, Geom. Funct.
  Anal., 11 (2001), pp.~465--588.

\bibitem{Green_Tao10}
{\sc B.~Green and T.~Tao}, {\em Linear equations in primes}, Ann. of Math. (2),
  171 (2010), pp.~1753--1850.

\bibitem{GT12-2}
{\sc B.~Green and T.~Tao}, {\em The {M}\"obius function is strongly orthogonal
  to nilsequences}, Ann. of Math. (2), 175 (2012), pp.~541--566.

\bibitem{GTZ12}
{\sc B.~Green, T.~Tao, and T.~Ziegler}, {\em An inverse theorem for the
  {G}owers {$U^{s+1}[N]$}-norm}, Ann. of Math. (2), 176 (2012), pp.~1231--1372.

\bibitem{Green61}
{\sc L.~W. Green}, {\em Spectra of nilflows}, Bull. Amer. Math. Soc., 67
  (1961), pp.~414--415.

\bibitem{HK02}
{\sc B.~Host and B.~Kra}, {\em An odd {F}urstenberg-{S}zemer\'edi theorem and
  quasi-affine systems}, J. Anal. Math., 86 (2002), pp.~183--220.

\bibitem{HK05-2}
\leavevmode\vrule height 2pt depth -1.6pt width 23pt, {\em Convergence of
  polynomial ergodic averages}, Israel J. Math., 149 (2005), pp.~1--19.
\newblock Probability in mathematics.

\bibitem{HK05}
\leavevmode\vrule height 2pt depth -1.6pt width 23pt, {\em Nonconventional
  ergodic averages and nilmanifolds}, Ann. of Math. (2), 161 (2005),
  pp.~397--488.

\bibitem{Host_Kra09}
\leavevmode\vrule height 2pt depth -1.6pt width 23pt, {\em Uniformity seminorms
  on {$\ell^\infty$} and applications}, J. Anal. Math., 108 (2009),
  pp.~219--276.

\bibitem{Leibman98}
{\sc A.~Leibman}, {\em Multiple recurrence theorem for measure preserving
  actions of a nilpotent group}, Geom. Funct. Anal., 8 (1998), pp.~853--931.

\bibitem{Leibman05-3}
\leavevmode\vrule height 2pt depth -1.6pt width 23pt, {\em Pointwise
  convergence of ergodic averages for polynomial actions of {${\Bbb Z}^d$} by
  translations on a nilmanifold}, Ergodic Theory Dynam. Systems, 25 (2005),
  pp.~215--225.

\bibitem{Leibman05}
\leavevmode\vrule height 2pt depth -1.6pt width 23pt, {\em Pointwise
  convergence of ergodic averages for polynomial sequences of translations on a
  nilmanifold}, Ergodic Theory Dynam. Systems, 25 (2005), pp.~201--213.

\bibitem{Leibman06}
\leavevmode\vrule height 2pt depth -1.6pt width 23pt, {\em Rational
  sub-nilmanifolds of a compact nilmanifold}, Ergodic Theory Dynam. Systems, 26
  (2006), pp.~787--798.

\bibitem{Leibman10}
\leavevmode\vrule height 2pt depth -1.6pt width 23pt, {\em Multiple polynomial
  correlation sequences and nilsequences}, Ergodic Theory Dynam. Systems, 30
  (2010), pp.~841--854.

\bibitem{Leibman10b}
\leavevmode\vrule height 2pt depth -1.6pt width 23pt, {\em Orbit of the
  diagonal in the power of a nilmanifold}, Trans. Amer. Math. Soc., 362 (2010),
  pp.~1619--1658.

\bibitem{Lesigne91}
{\sc E.~Lesigne}, {\em Sur une nil-vari\'et\'e, les parties minimales
  associ\'ees \`a une translation sont uniquement ergodiques}, Ergodic Theory
  Dynam. Systems, 11 (1991), pp.~379--391.

\bibitem{Parry69}
{\sc W.~Parry}, {\em Ergodic properties of affine transformations and flows on
  nilmanifolds.}, Amer. J. Math., 91 (1969), pp.~757--771.

\bibitem{Parry70}
\leavevmode\vrule height 2pt depth -1.6pt width 23pt, {\em Dynamical systems on
  nilmanifolds}, Bull. London Math. Soc., 2 (1970), pp.~37--40.

\bibitem{Raghunathan72}
{\sc M.~S. Raghunathan}, {\em Discrete subgroups of {L}ie groups},
  Springer-Verlag, New York-Heidelberg, 1972.
\newblock Ergebnisse der Mathematik und ihrer Grenzgebiete, Band 68.

\bibitem{Rauzy76}
{\sc G.~Rauzy}, {\em Propri\'et\'es statistiques de suites arithm\'etiques},
  Presses Universitaires de France, Paris, 1976.
\newblock Le Math{\'e}maticien, No. 15, Collection SUP.

\bibitem{Ribenboim96}
{\sc P.~Ribenboim}, {\em The new book of prime number records},
  Springer-Verlag, New York, 1996.

\bibitem{Sun15}
{\sc W.~Sun}, {\em Multiple recurrence and convergence for certain averages
  along shifted primes}, Ergodic Theory Dynam. Systems, 35 (2015),
  pp.~1592--1609.

\bibitem{Vinogradov47}
{\sc I.~M. Vinogradov}, {\em The method of trigonometrical sums in the theory
  of numbers}, Interscience Publishers, London and New York., 1947.
\newblock Translated, revised and annotated by K. F. Roth and Anne Davenport.

\bibitem{Ziegler05}
{\sc T.~Ziegler}, {\em A non-conventional ergodic theorem for a nilsystem},
  Ergodic Theory Dynam. Systems, 25 (2005), pp.~1357--1370.

\bibitem{Ziegler07}
{\sc T.~Ziegler}, {\em Universal characteristic factors and {F}urstenberg
  averages}, J. Amer. Math. Soc., 20 (2007), pp.~53--97 (electronic).

\end{thebibliography}


\begin{thebibliography}{1}

\bibitem{BHK05}
{\sc V.~Bergelson, B.~Host, and B.~Kra}, {\em Multiple recurrence and
  nilsequences}, Invent. Math., 160 (2005), pp.~261--303.
\newblock With an appendix by Imre Ruzsa.

\bibitem{Leibman06}
{\sc A.~Leibman}, {\em Rational sub-nilmanifolds of a compact nilmanifold},
  Ergodic Theory Dynam. Systems, 26 (2006), pp.~787--798.

\bibitem{MR19}
{\sc J.~Moreira and F.~K. Richter}, {\em A spectral refinement of the
  {B}ergelson-{H}ost-{K}ra decomposition and new multiple ergodic theorems},
  Ergodic Theory and Dynamical Systems, 39 (2019), pp.~1042--1070.

\end{thebibliography}
% ================================================================

% ==================================================END DOCUMENT
\end{document}